\documentclass{amsart}
\usepackage{amssymb,amscd}   
   
\begin{document}   
    
%
%
\setlength{\marginparwidth}{1.12in}   
\newcommand{\mar}[1]{{\marginpar{\textsf{#1}}}}   
\newcommand\datver[1]{\def\datverp%
 {\par\boxed{\boxed{\text{Version: #1; Run: \today}}}}}   
\datver{4.3; Revised: 7/26/2000}   
   
%
%
\newcommand\Dr{{\not \!\!D}} 
\newcommand\Dir{{\not \!\!D}} 
\newcommand{\wt}[1]{\widetilde{#1}} { }\vspace*{-4mm}  
\newcommand{\Prod}{\prod}   
\newcommand{\Cal}{\mathcal}   
\newcommand\Ran{\operatorname{Im}}   
\newcommand\abs{\operatorname{abs}}   
\newcommand\rel{\operatorname{rel}}   
\newcommand\inv{\operatorname{inv}}   
\newcommand\topo{\operatorname{top}}   
\newcommand\opp{\operatorname{op}}   
\newcommand\mfk{\mathfrak}   
\newcommand\coker{\operatorname{coker}}   
\newcommand\hotimes{\hat \otimes}   
\newcommand\ind{\operatorname{ind}}   
\newcommand{\id}{\operatorname{id}}   
\newcommand\End{\operatorname{End}}   
\newcommand\per{\operatorname{per}}   
\newcommand\pa{\partial}   
\newcommand\sign{\operatorname{sign}}   
\newcommand\supp{\operatorname{supp}}   
\newcommand{\comp}{\operatorname{comp}}   
\newcommand{\ck}{{\mathcal K}}   
\newcommand{\cb}{{\mathcal B}}   
\newcommand{\cD}{\mathcal{D}}   
\newcommand{\cH}{\mathcal{H}}   
\newcommand{\cV}{\mathcal{V}}   
\newcommand{\cHmu}{\cH_{-\infty}}   
\newcommand{\cHu}{\cH_{\infty}}   
\newcommand\cy{\mathcal{C}^\infty}   
\newcommand\CI{\mathcal{C}^\infty}   
\newcommand\CO{\mathcal{C}_0}  
\newcommand{\cunc}{\CI_{c}}   
\newcommand\lra{\longrightarrow}   
\newcommand\vlra{-\!\!\!-\!\!\!-\!\!\!\!\longrightarrow}   
\newcommand\bS{{}^b\kern-1pt S}   
\newcommand\bT{{}^b\kern-1pt T}   
\newcommand\Hom{\operatorname{Hom}}   
\newcommand{\adb}{\operatorname{ad}}   
\newcommand{\norm}[1]{\| #1 \|}   
\newcommand{\rs}{\widetilde{\varrho}}    
   
\newcommand\alg[1]{\mathfrak{A}(#1)}    
\newcommand\qalg[1]{\mathfrak{B}(#1)}    
\newcommand\ralg[1]{\mathfrak{A}_r(#1)}    
\newcommand\rqalg[2]{\mathfrak{B}_r(#1)}   
\newcommand\ideal[1]{C^*(#1)}    
\newcommand\rideal[1]{C^*_r(#1)}    
\newcommand\qideal[2]{\mathfrak{R}_{#1}(#2)}    
\newcommand\In{\operatorname{In}}

\newcommand\TR{\operatorname{T}}   
\newcommand\ha{\frac12}    
\newcommand\cal{\mathcal}   
\newcommand\END{\operatorname{END}}   
\newcommand\ENDG{\END_{\GR}(E)}   
\newcommand{\eTM}{{}^{e}TM}   
\newcommand\CC{\mathbb C}   
\newcommand\NN{\mathbb N}   
\newcommand{\nzn}{\NN_{0}}   
\newcommand\RR{\mathbb R}   
\newcommand\ZZ{\mathbb Z}   
\newcommand\ci{${\mathcal C}^{\infty}$}   
\newcommand\CIc{{\mathcal C}^{\infty}_{\text{c}}}   
\newcommand\hden{{\Omega^{\lambda}_d}}   
\newcommand\VD{{\mathcal D}}   
\newcommand{\VeM}{\mathcal{V}_{e}(M)}   
\newcommand\mhden{{\Omega^{-1/2}_d}}   
\newcommand\ehden{r^*(E)\otimes {\Omega^{\lambda}_d}}   
   
   
\newcommand{\Cat}{\mathcal C}   
\newcommand{\Cliff}{{\rm Cliff}}   
\newcommand{\Diff}[1]{{\rm Diff}(#1)}   
\newcommand{\Gr}[1]{{\mathcal G}^{(#1)}}   
\newcommand{\GR}{\mathcal G}   
\newcommand{\LGR}{\mathcal L}   
\newcommand{\BB}{\mathbb{B}}   
\newcommand{\GG}{\mathcal G}   
\newcommand{\OA}{\mathcal O}   
\newcommand{\lgg}{\mathfrak g}   
\newcommand{\tPS}[1]{\Psi^{#1}(\GR)}   
\newcommand{\PS}[1]{\Psi^{#1}(\GR;E)}   
\newcommand{\PSF}[1]{\Psi^{#1}(\GR;F)}   
\newcommand{\AL}{{\mathcal A}(\GR)}   
\newcommand{\FAM}{P=(P_x)\,,x \in M,}   
\newcommand\symb[2]{{\mathcal S}^{#1}(#2)}   
\newcommand{\loc}{\operatorname{loc}}   
\newcommand{\cl}{\operatorname{cl}}    
\newcommand{\mI}{\mathfrak I}     
\newcommand{\MGR}{\mathcal M}

   
%
%
\newcommand\Mand{\text{ and }}   
\newcommand\Mandset{\text{ and set}}   
\newcommand\Mas{\text{ as }}   
\newcommand\Mat{\text{ at }}   
\newcommand\Mfor{\text{ for }}   
\newcommand\Mif{\text{ if }}   
\newcommand\Min{\text{ in }}   
\newcommand\Mon{\text{ on }}   
\newcommand\Motwi{\text{ otherwise }}   
\newcommand\Mwith{\text{ with }}   
\newcommand\Mwhere{\text{ where }}    
\newcommand\ie{{\em i.e., }}    
%
%
\newtheorem{theorem}{Theorem}   
\newtheorem{proposition}{Proposition}   
\newtheorem{corollary}{Corollary}   
\newtheorem{lemma}{Lemma}   
\newtheorem{definition}{Definition}   
\newtheorem{notation}{Notations}   
\theoremstyle{remark}   
\newtheorem{remark}{Remark}   
\newtheorem{example}{Example}

%
%
   
\author[R. Lauter]{Robert Lauter} \address{Universit\"at   
       Mainz. Fachbereich 17-Mathematik, D-55099 Mainz, Germany   
       \newline   
       {\it currently}:   
        Universit\"at M\"unster, SFB 478,   
                Hittorfstra{\ss}e 27,   
                D-48149 M\"unster, Germany}   
       \email{lauter@mathematik.uni-mainz.de}

\author[V. Nistor]{Victor Nistor} \address{Pennsylvania State   
       University, Math. Dept., University Park, PA 16802}   
       \email{nistor@math.psu.edu}    
   
\thanks{Lauter was partly supported by a scholarship of the   
        German Academic Exchange Service (DAAD) within  
        the {\em Hochschulsonderprogramm III von Bund und  
        L\"andern}, and the   
        Sonderforschungsbereich 478  
        {\em Geometrische Strukturen in der Mathematik}  
        at the University of M\"unster.    
      Nistor was partially supported an  NSF Young   
      Investigator Award DMS-9457859  and NSF Grant DMS-9971951. {\bf   
        http:{\scriptsize//}www.math.psu.edu{\scriptsize/}nistor{\scriptsize/}}   
       .}   
    
\dedicatory\datverp   
\begin{abstract}   
The first five sections of this paper are a survey of algebras of 
pseudodifferential operators on groupoids. We thus review 
differentiable groupoids, the definition of pseudodifferential 
operators on groupoids, and some of their properties. We use then this 
background material to establish a few new results on these algebras 
that are useful for the analysis of geometric operators on non-compact 
manifolds and singular spaces. The first step is to establish that the 
geometric operators on groupoids are in our algebras. This then leads 
to criteria for Fredholmness for geometric operators on suitable 
non-compact manifolds, as well as to an inductive procedure to study 
their essential spectrum.  As an application, we answer a question of 
Melrose on the essential spectrum of the Laplace operator on manifolds 
with multi-cylindrical ends. 
\end{abstract}    
   
\title[Geometric operators]{Analysis of geometric operators on open    
manifolds:\ A groupoid approach}    
\maketitle   
\tableofcontents

\def\c{\cite}   
\def\fr{\frac}   
\def\ub{\underbar}   
\def\O{\Cal O}   
\def\F{\Cal F}   
\def\differ{\text{differentiable} }   
\def\tPSeudo{\text{pseudodifferential} }   
\def\supp{\text{supp} }   
\def\inn{{\mathcal R}}    
\def\prop{\text{prop}}    
\def\frag{\frak{G}}   
\def\simd{\tilde{d}}   
\def\simf{\tilde{\F}}   
\def\simo{\tilde{\O}}   
\def\simr{\tilde{r}}   
\def\simp{\tilde{p}}   
\def\simmu{\tilde{\mu}}

\section*{Introduction}

The first half of this paper is a survey of the results from  
\cite{LMN,NistorInt,NistorIndFam} and \cite{NWX}. However, there are  
some new results, including a determination of the essential spectrum  
of the Laplace operator on complete manifolds with multi-cylindrical  
ends.  This was formulated as a question in \cite{MelroseScattering}  
(Conjecture 7.1).  
   
Let us now discuss the contents of this paper. As we mentioned above,  
the first five sections of this paper are mostly a survey of results  
on pseudodifferential operators on groupoids. In Section  
\ref{Sec.prelim}, we review some definitions involving manifolds with  
corners and we introduce groupoids. We also define in this section the  
class of groupoids we are interested in, namely ``differentiable  
groupoids'' (Definition \ref{Almost.differentiable}), and we recall  
the definition of the Lie algebroid associated to a differentiable  
groupoid. (Our ``differentiable groupoids'' should more properly be  
called ``Lie groupoids.'' However, this term was already used for some  
specific classes of differentiable groupoids.) The construction which  
associates to a differentiable groupoid $\GR$ its Lie algebroid  
$A(\GR)$ is a generalization of the construction which associates to a  
Lie group its Lie algebra.  This construction is made possible for  
differentiable groupoids by the fact that the fibers $\GR_x :=  
d^{-1}(x)$ of the domain map $d : \GR \to M$ consist of smooth  
manifolds without corners (in this paper, a ``smooth manifold'' will  
always mean a ``smooth manifold without corners,'' and a  
$C^\infty$-manifold with corners will be called simply, a  
``manifold''). It is convenient to think of $\GR$ as a set of arrows  
between various points, called units, which can be composed according  
to some definite rules. If $g\in \GR$ is such an ``arrow,'' then  
$d(g)$ is its domain and $r(g)$ is its range, so $\GR_x$ is the set of  
all arrows starting at (or with domain) $x$.  
   
Section \ref{Sec.Examples.I} contains several examples of  
differentiable groupoids.  In Section \ref{Sec.POOG}, we introduce  
pseudodifferential operators on groupoids. A pseudodifferential  
operator $P$ on the differentiable groupoid $\GR$ is actually a family  
$P=(P_x)$ of ordinary pseudodifferential operators $P_x$ on the smooth  
manifolds without corners $\GR_x$. This family is required to be  
invariant with respect to the natural action of $\GR$ by right  
translations and to be differentiable in a natural sense. Because the  
family $P$ acts on $\CIc(\GR)$, we can regard $P$ as an operator on  
this space.  To get the functoriality and composition properties  
right, we also assume that the convolution kernel $\kappa_P$ of $P$ is  
compactly supported ($\kappa_P$ has been called ``the reduced kernel of  
$P$'' in \cite{NWX}). We denote by $\tPS{m}$ the space of such order  
$m$ (families of) pseudodifferential operators on $\GR$. Many  
properties of the usual pseudodifferential operators on smooth  
manifolds extend to the operators in $\tPS{m}$, most notably, we get  
the existence of the principal symbol map. At the end of the section,  
we indicate how to treat operators acting between sections of two  
vector bundles. Section \ref{Sec.BR} deals with the necessary facts  
about the actions of $\tPS{\infty}$ on various classes of functions.  
It is convenient to present this from the point of view of  
representation theory (after all, this is the representation theory of  
the Lie group $G$ if $\GR = G$).  It is in fact enough to study  
representations of $\tPS{-\infty}$, and this includes the study of  
boundedness properties of various representations: to check that a  
representation of $\tPS{0}$ consists of bounded operators, it is  
sufficient to check that its restriction to $\tPS{-\infty}$ is  
bounded.  This generalizes the usual boundedness theorems   
pseudodifferential operators of order $0$.  
   
Motivated by a question of Connes, Monthubert also was lead to define  
pseudodifferential operators on groupoids in \cite{Monthubert},  
independently. Actually, Connes had defined algebras of  
pseudodifferential operators on the graphs of $C^{0,\infty}$  
foliations in \cite{ConnesF}, before. The definitions of  
pseudodifferential operators on groupoids in \cite{Monthubert} and  
\cite{NWX} is essentially the same as Connes', but different because  
we take into account the differentiability in the transverse  
direction, too, and, most importantly, we allow non-regular Lie  
algebroids.  An approach to operators on singular spaces which is  
similar in philosophy to ours was outlined in  
\cite{MelroseScattering}.  There Melrose considers operators whose  
kernels are defined on compact manifolds with corners and have  
structural maps that make them ``semi-groupoids.''  The results of  
\cite{NWX} were first presented in July 1996 at the joint SIAM-AMS-MAA  
Meeting on Quantization in Mount Holyoke.  
   
In addition to the survey of the results from   
\cite{LMN,NistorInt,NistorIndFam} and \cite{NWX}, the first five   
sections of the paper also include many examples of groupoids together   
with the description of the resulting algebras of pseudodifferential   
operators. Actually, Sections \ref{Sec.Examples.I} and   
\ref{Sec.Examples.II} are devoted exclusively to examples, with the hope   
that this will make the general theory easier to apply. Whenever it   
was relevant, we have also compared our construction to the classical   
constructions.   
   
Let us now quickly describe the contents of each of the remaining five  
sections of this paper. In Section \ref{Sec.GO} we show that the  
geometric differential operators acting on the fibers of the domain  
map $d:\GR \to M$ belong to our algebras $\PS{\infty}$, for suitable  
$E$. To define these geometric operators -- except the de Rham  
operator -- we need a metric on $\GR_x$, and this will come from a  
metric on $A(\GR)$, the Lie algebroid of $\GR$. Then we need to  
establish the existence of right invariant connections with the  
properties necessary to define the geometric operators we are  
interested in. It needs to be established, for example, that  
compatible connections exist on Clifford modules, and this is not  
obvious in the groupoid case. Section \ref{Sec.Sobolev} establishes  
the technical facts needed to define Sobolev spaces in our  
setting. The reader should find the statements in that section easy to  
understand and believe (although, unfortunately, not the same thing  
can be said about their proofs).  
   
Beginning with Section \ref{Sec.O.M}, we begin to work with groupoids   
$\GR$ of a special kind, which model operators on suitable non-compact   
Riemannian manifolds $(M_0,g)$. The set of units $M$ of these groupoids is a   
compactification of $M_0$ to a manifold with corners, and hence we can   
think of $\GR$ as modeling the behavior at $\infty$ of $M_0$ (this    
approach was inspired in part by Melrose's geometric scattering theory   
program). More precisely, we assume that $M_0$   
is an open invariant subset of $M$ with the property that $\GR_{M_0}$   
is the pair groupoid, that is, that for each $x \in M_0$ there exists   
exactly one arrow to any other point $y \in M_0$, and there exists no   
arrow to any point not in $M_0$ (see Example \ref{ex.pair}). Then the   
restriction of $A(\GR)$ to $M_0$ identifies canonically with the   
tangent bundle to $M_0$. For such Riemannian manifolds, the geometric   
operators on $M_0$ can be recovered from the geometric operators on   
$\GR$. We use these results in Section \ref{Sec.Examples.III} to study   
the  Hodge-Laplace operators on suitable non-compact Riemannian   
manifolds, using the general spectral properties discussed in Section   
\ref{Sec.SP}. We obtain, in particular, criteria for certain   
pseudodifferential operators on $M_0$ to be Fredholm or compact,   
similar to the well-known criteria for $b$-pseudodifferential   
operators to be Fredholm or compact \cite{Melrose-Nistor1,mepi92}

\subsubsection*{Acknowledgments:}    
We have benefited from discussions and suggestions from several 
people: A. Connes, R. Melrose, B. Monthubert, and G. Skandalis, and we 
would like to thank them. We are grateful to R. Melrose who has 
instructed us on the $b$-calculus. The first author would like to 
thank R. Melrose, J. Cuntz, and the SFB 478 {\em Geometrische 
Strukturen in der Mathematik}, for their warm hospitality and useful 
discussions at MIT and, respectively, at the University of M\"unster, 
where part of this work was done.

\section{Manifolds with corners and groupoids\label{Sec.prelim}}

In the following, we shall consider manifolds with corners.  Brief 
introductions into the analysis on manifolds with corners can be found 
for instance in \cite{mepi92}, and will be discussed in more detail 
in a forthcoming book of Melrose. We begin this section with a short 
discussion of the relevant definitions concerning manifolds with 
corners. Then we introduce groupoids and the class of {\em 
differentiable groupoids}. The reader can find more information on 
groupoids in \cite{connes, Renault1}.  
   
By a {\em manifold} we shall always understand a differentiable  
manifold possibly {\em with corners}, by a {\em smooth manifold} we  
shall always mean a differentiable manifold {\em without corners}. By  
definition, every point $p$ in a manifold with corners $M$ has a  
coordinate neighborhood diffeomorphic to $[0,\infty)^k \times  
\RR^{n-k}$ such that the transition functions are smooth (including on  
the boundary). Moreover, we assume that each boundary hyperface $H$ of  
$M$ is an embedded submanifold and has a defining function, that is,  
that there exists a smooth function $x \ge 0$ on $M$ such that  
$$   
        H = \{ x = 0 \} \text{ and } dx \not = 0 \text{ on }H\,.$$   
It follows that if $H_1,\ldots,H_k$ are boundary   
hyperfaces of $M$ with defining functions   
$x_1,\ldots,x_k$, then the differentials   
$dx_{1}, \ldots, dx_k$ are linearly independent at   
the intersection $H_1 \cap \ldots \cap H_k$.

\begin{definition}\label{def.1}\   
A \underline{submersion} $f : M \to N$, between two manifolds with corners   
$M$ and $N$, is a differentiable map such that $df$ is surjective   
at all points and $df(v)$ is an inward pointing tangent vector of $N$   
if, and only if, $v$ is an inward pointing vector $M$.    
\end{definition}

It is not difficult to prove the following lemma.

\begin{lemma}\ Let $f : M \to N$ be a submersion of manifolds with corners,    
$y \in N$ a point belonging to the interior of a face of    
codimension $k$ and $x_1,x_2,\ldots,x_k$ be the defining functions of   
the hyperfaces containing $y$. Then $x_1\circ f, x_2\circ f, \ldots,   
x_k \circ f$ are defining functions for $k$ distinct hyperfaces of $M$.   
Let $p \in M$ be such that $f(p) = y$, then all   
hyperfaces of $M$ containing $p$ are obtained in this way.   
\end{lemma}

\begin{proof}\ Let $p$ and $y$ be as in the statement and   
$z_j := x_j \circ f$. Because $df$ is surjective, it follows that $dz_j$   
are linearly independent at $p$. Each function $z_j$ is the   
defining function of a hyperface because $f$ must map faces of   
codimension $k$ to faces of codimension $k$.   
\end{proof}

For any submersion $f$ as above, it follows    
that the fibers $f^{-1}(y)$ of $f$ are smooth manifolds   
(that is {\em without}   
corners). We can see this as follows. Because this property of $f$    
is local, we can fix $y$ and replace $M$ and $N$ with some small   
open neighborhoods   
of $y$ and $f(y)$, respectively. By decreasing these neighborhoods,   
we can also assume that they are diffeomorphic to one of the model open sets   
$[0,1)^k \times \RR^{n-k}$. Then we extend $M$ and $N$ to smooth   
open manifolds without corners, denoted $\tilde M$ and $\tilde N$, and   
$f$ to a smooth function $\tilde f : \tilde M \to \tilde N$. By decreasing   
$M$ and $N$ again, if necessary, we can assume that $\tilde f$ is still   
a submersion, this time in the classical sense, because $\tilde M$ and   
$\tilde N$ are smooth manifolds. This gives then that   
$\tilde f^{-1}(y)$ is a smooth submanifold of $\tilde{M}$.   
Since $f^{-1}(y) = \tilde f^{-1}(y) \cap M$, it is enough to check that    
this intersection coincides with a component of $\tilde f^{-1}(y)$.   
Let $x$ be a defining function of $N$ in $\tilde N$ (that is,   
the defining function of a hyperface of $N$). Then $x \circ f$ is   
a defining function of some hyperface of $M$. By counting the hyperfaces of   
$M$ and $N$, we see that we get in this way all defining functions of $M$.   
Since they all vanish on $y$, $f^{-1}(y) \cap \{x \circ f \ge 0 \} =   
f^{-1}(y)$. This completes the argument.

The concept of a ``submanifold'' means the following in our setting.

\begin{definition}\label{def.2}\   
A \underline{submanifold} (or \underline{submanifold with corners})    
$N$ of a manifold with corners $M$   
is a submanifold $N \subset M$ such that $N$ is a manifold with corners and   
each hyperface $F$ of $N$ is a connected component of a set   
of the form $F' \cap N$, where $F'$   
is a hyperface of $M$ intersecting $N$ transversally.   
\end{definition}

We shall need groupoids endowed with various   
structures. (\cite{Renault1} is a general reference for some of what   
follows.)  We recall first that a {\em small category} is a category whose   
class of morphisms is a set.  The class of objects of a small category   
is then a set as well. Here is now a quick definition of groupoids.

\begin{definition}\label{def.3}   
A \underline{groupoid} is a small category $\GR$ all of whose    
morphisms are invertible.    
\end{definition}

Let us make this definition more explicit. A groupoid $\GR$ is   
a pair $(\Gr0,\Gr1)$ of sets together  with   
structural morphisms $d,r,\mu,u$, and $\iota$.   
Here the first set, $\Gr0$, represents the   
objects (or units) of the groupoid and the second set, $\Gr1$, represents    
the set of morphisms of $\GR$. Usually, we shall denote the space of units   
of $\GR$ by $M$ and we shall identify $\GR$ with $\Gr1$. Each object of    
$\GR$ can be identified with a morphism of $\GR$, the identity morphism of   
that object, which leads to an injective map $u : M := \Gr0 \to \GR$,   
used to identify $M$ with a subset of $\GR$.   
Each morphism $g \in \GR$ has a ``domain'' and a ``range.''   
We shall denote by $d(g)$     
the {\em domain} of  $g$ and by $r(g)$ the {\em range} of $g$.   
We thus obtain functions    
\begin{equation}   
        d,r:  \GR \longrightarrow  M:=\Gr0.   
\end{equation}   
The  multiplication (or composition) $\mu(g,h)=gh$ of two morphisms    
$g, h \in \GR$ is   
not always defined; it is defined precisely when $d(g) = r(h)$.     
The inverse of a morphism $g$   
is denoted by $g^{-1}=\iota(g)$.   
   
A  groupoid $\GR$ is completely determined by    
the spaces $\Gr0$ and $\Gr1$ and the structural maps $d,r,\mu,u,\iota$.   
We sometimes write $\GR=(\Gr0,\Gr1,d,r,\mu,u,\iota)$. The structural   
maps satisfy the following properties:   
\\   
\begin{itemize}    
\item  $r(gh)=r(g)$ and $d(gh)=d(h)$,   
       for any pair $g,h$ with $d(g) = r(h)$;\\   
\item  The partially defined multiplication $\mu$ is associative;\\   
\item  $d(u(x))=r(u(x))=x$, $\forall x\in \GR^{(0)}$,\ $u(r(g))g=g$,    
$gu(d(g))=g$,  $\forall g\in \GR^{(1)}$,\ and $u:\Gr0 \to \Gr1$   
is one-to-one;\\   
\item  $r(g^{-1})=d(g)$,\ $d(g^{-1})=r(g)$,\ $gg^{-1}=u(r(g))$,\ and    
$g^{-1}g=u(d(g))$. \\  
\end{itemize}

\begin{definition}\label{Almost.differentiable}   
A \underline{differentiable groupoid} is a groupoid    
$$   
        \GR=(\Gr0,\Gr1,d,r,\mu,u,\iota)   
$$   
such that $M:=\Gr 0$ and $\Gr1$ are manifolds with corners,    
the structural maps $d,r,\mu,u,$ and $\iota$ are differentiable,    
the domain map $d$ is a submersion, and all the spaces    
$M$ and $\GR_x := d^{-1}(x)$, $x\in M$,  are Hausdorff.    
\end{definition}

Note that we {\em do not} require $\Gr1$ to be Hausdorff. We actually   
need this generality to treat algebras associated to foliations   
and other geometric structures in our setting.    
   
We also observe that $\iota$ is a diffeomorphism, and hence $d$ is a   
submersion if, and only if, $r=d\circ \iota$ is a submersion.     
The reason for requiring $d$ to be a submersion is that then   
each fiber $\GR_x=d^{-1}(x) \subset   
\Gr1$ is a smooth manifold without corners.   
   
We now recall the definition of the ``Lie algebroid'' of a    
differentiable   
groupoid. Lie algebroids are for differentiable groupoids   
what Lie algebras are for Lie groups.

\begin{definition}\label{Lie.Algebroid} A  \underline{Lie algebroid} $A$     
over a manifold $M$ is a vector bundle $A$ over $M$, together with a Lie   
algebra structure on the space $\Gamma(A)$ of smooth sections of $A$   
and  a bundle map $\varrho: A \rightarrow TM$, extended to a map between    
sections of these bundles,  such that   
   
(i) $\varrho([X,Y])=[\varrho(X),\varrho(Y)]$, and   
   
(ii) $[X, fY] = f[X,Y] + (\varrho(X) f)Y$,   
   
\noindent for any smooth sections $X$ and $Y$ of $A$ and any smooth function   
$f$ on $M$.   
The map $\varrho$ is called the {\em anchor}.   
\end{definition}

Note that we allow the base $M$ in the definition above to be a    
manifold with corners.    
   
The Lie algebroid associated to a differentiable groupoid $\GR$   
is defined as follows \cite{Mackenzie1,Pradines2}. The vertical tangent    
bundle (along the fibers of $d$)   
of a differentiable groupoid $\GR$ is, as usual,    
\begin{equation}   
        T_{vert} \GR = \ker d_*    
        = \bigcup_{x\in M} T \GR_x \subset T\GR.   
\end{equation}   
Then $A(\GR) := T_{vert} \GR\big |_{M}$, the restriction of   
the $d$-vertical tangent bundle to the set of units, $M$, determines   
$A(\GR)$ as a vector bundle.   
   
We now construct the bracket defining the Lie algebra structure on $A(\GR)$.   
The right translation by an arrow $g \in \GR$ defines   
a diffeomorphism    
$$   
        R_g:\GR_{r(g)}\ni g' \longmapsto g'g \in \GR_{d(g)}.   
$$     
A vector field $X$ on $\GR$ is called $d$-vertical if $d_*(X(g)) = 0$ for    
all $g$. The $d$-vertical vector fields are precisely the vector fields   
on $\GR$ that can be restricted to the submanifolds $\GR_x$. It makes sense   
then to consider right--invariant vector fields on $\GR$.    
It is not difficult to see that the sections of $A(\GR)$ are   
in one-to-one correspondence with $d$--vertical, right--invariant   
vector fields on $\GR$.    
    
The Lie bracket $[X,Y]$   
of two $d$--vertical right--invariant vector fields $X$ and $Y$ is   
also $d$--vertical and right--invariant, and hence the Lie bracket   
induces a Lie algebra structure on the sections of $A(\GR)$. To define   
the action of the sections of $A(\GR)$ on functions on $M$, observe   
that the right invariance property makes sense also for functions on   
$\GR$, and that $\CI(M)$ may be identified with the subspace of smooth   
right--invariant functions on $\GR$, because $r$ is a submersion.    
If $X$ is a right--invariant   
vector field on $\GR$ and $f$ is a right--invariant function on $\GR$,   
then $X(f)$ will still be a right invariant function. This identifies   
the action of $\Gamma(A(\GR))$ on $\CI(M)$.   
   
We denote by $T_{vert}^* \GR$ the dual of $T_{vert} \GR$, and by $A^*(\GR)$    
the dual of $A(\GR)$.     
Later on, we shall need the bundle $\hden$ of $\lambda$-densities along    
the fibers of $d$. It is defined as follows. If the   
fibers of $d$ have dimension $n$, then    
$$   
        \hden :=|\Lambda^n T_{vert}^*\GR|^\lambda =  
        \cup_x\Omega^\lambda(\GR_x).   
$$    
By invariance, these bundles can be obtained as pull-backs of bundles   
on $M$. For example $T_{vert} \GR = r^*(A(\GR))$ and $\hden   
=r^*({\mathcal D}^\lambda)$, where ${\mathcal D}^\lambda$ denotes   
$\hden\vert_{M}$.  If $E$ is a (smooth complex) vector bundle on the   
set of units $M$ of $\GR$, then the pull-back bundle $r^*(E)$ on $\GR$   
will have right invariant connections obtained as follows. A   
connection $\nabla$ on $E$ lifts to a connection on $r^*(E)$. Its   
restriction to any fiber $\GR_{x}$ defines a linear connection in the   
usual sense, which is denoted by $\nabla_{x}$. It is easy to see that   
these connections are right invariant in the sense that   
\begin{equation}   
\label{eq.invariant.conn}   
        R_{g}^*\nabla_{x}=\nabla_{y}, \  \ \ \forall g\in \GR   
        \mbox{ such that } r(g)=x \mbox{ and } d(g)=y.   
\end{equation}   
The bundles considered above will thus have invariant connections.   
   
We observe that in all considerations above, we first use the smooth   
structure on each $\GR_x$ to define the geometric quantities we are   
interested in:\ $X(f)$, $[X,Y]$, and so on. We do need then, however,   
to check that these quantities define global objects on the possibly   
non-Hausdorff manifold $\GR$, more precisely, we need the defined   
objects to be smooth on $\GR$, not just on every $\GR_x$. All these   
global smoothness conditions can be checked on {\em smooth} functions   
on $\GR$, as long as we correctly define this concept. For   
non-Hausdorff manifolds, the correct choice is the one considered by
Crainic and Moerdijk in \cite{cramo00}, more precisely, we consider first   
the spaces   
$V=\oplus_\alpha \CIc(U_\alpha)$ and $W = \oplus_{\alpha,\beta}   
\CIc(U_\alpha \cap U_\beta)$, where $U_\alpha \subset \GR$ are the   
domains of coordinate charts. Then there is a natural map $\delta:W \to V$,   
which is the direct sum of the maps $(j,-j) : \CIc(U_\alpha \cap   
U_\beta) \to \CIc(U_\alpha) \oplus \CIc(U_\beta)$, with $j$ the   
natural inclusion, and we define  
\begin{equation}\label{eq.Crainic}   
        \CIc(\GR) = V/\delta(W).   
\end{equation}   
A function $f$ on $\GR$ is {\em smooth} if, and only if,    
$\phi f \in \CIc(U_{\alpha})$ for all   
$\phi \in \CIc(U_\alpha)$.    
   
If $A \to M$ is a given Lie algebroid and $\GR$ is a differentiable  
groupoid whose Lie algebroid is isomorphic to $A$, then we say that  
$\GR$ {\em integrates} $A$. Not every Lie algebroid can be integrated  
(see \cite{al-mo:suites} for an example). There are many simple-minded  
Lie algebroids for which the standard integration procedures of  
\cite{NistorInt} lead to non-Hausdorff groupoids. For many Lie  
algebroids (including the tangent bundles to foliations) there are no  
Hausdorff groupoids integrating them. The integration of Lie algebroids  
was also studied by Moerdijk in \cite{moerdijk}.

\section{Examples I: Groupoids\label{Sec.Examples.I}}

We now discuss examples of differentiable groupoids. The   
examples included in this section are either basic examples or theoretical   
examples predicted by the general theory of groupoids. The   
examples that we shall use to study geometric operators on open   
manifolds will be included in Section \ref{Sec.Examples.II}.    
For each example considered in this section we also   
describe the corresponding Lie algebroid.    
   
As in the previous section, we shall identify    
$\GR$ with its set of arrows: \ $\GR = \Gr1$, and we shall denote by   
$M$ or by $\Gr0$ the set of units of $\GR$.    

\begin{example}\label{ex.manifold}\    
{\em Manifolds with corners:\ }\ If $M$ is a manifold with   
corners, then $M$ is naturally a differentiable groupoid with no   
arrows other than the units, \ie we have $\GR = \Gr1 = \Gr0=M $, $d=r=id$.    
   
Then $A(\GR) = 0$, the zero bundle on $M$.   
\end{example}

\begin{example}\label{ex.Lie}\ {\em Lie groups:\ }\ Every differentiable    
groupoid  $\GR$   
with space of units consisting of just one point, $M = \{e\}$, is necessarily   
a Lie group. Conversely, every Lie group $G$ can be regarded as an    
differentiable   
groupoid $\GR = G$ with exactly one unit $M = \{e\}$, the unit of that group.    
   
In this case $A(\GR) = Lie (G)$, the Lie algebra of $G$.   
\end{example}

The above two sets of examples are in a certain way the two extreme  
examples of differentiable groupoids. The first set of examples  
consists of groupoids each of which has the least set of arrows among  
all groupoids with the same set of units. The second set of examples  
consists of the differentiable groupoids with the least set of units  
among all non-empty groupoids. The fact that manifolds and Lie groups  
are particular cases of groupoids makes them a favorite object of  
study in non-commutative geometry.

The following example plays a special role in the theory pseudodifferential   
operators on groupoids.

\begin{example}\label{ex.pair}\ {\em The pair groupoid:\ }\   
Let  $M$ be a smooth manifold (without corners)   
and let   
$$   
        \GR  = M \times M \quad \text{and} \quad \Gr0 = M,   
$$   
with structural morphisms    
$d(x,y)=y$, $r(x,y)=x$, $(x,y)(y,z)=(x,z)$, $u(x) = (x,x)$    
and $\iota(x,y) = (y,x)$.   
Then $\GR$ is a differentiable groupoid, called    
{\em the pair groupoid}.    
   
Denote the pair groupoid with units $M$ by $M \times M$.  Then $A(M  
\times M) = TM$, the tangent bundle to $M$.  
\end{example}

A variant of the above example is the following.

\begin{example}\label{ex.fibered.pair}\    
{\em The fibered pair groupoid:\ }\ Let $p : Y \to B$ be a submersion  
of manifolds with corners (see Definition \ref{def.1} of the previous  
section).  The fibered pair groupoid $\GR$ is obtained as  
$$   
        \GR = Y \times_B Y := \{(x,y), p(x)=p(y), x,y \in Y \},   
$$    
with the operations induced from the pair groupoid $Y \times Y$. Its space   
of units is $Y$.   
   
The Lie algebroid of $Y \times_B Y$ is $A(\GR) = T_{vert}Y$,   
the kernel of $TY \to TB$, or, in other words, the vertical tangent bundle    
to the submersion $p : Y \to B$.   
\end{example}

\begin{example}\ {\em Products:\ }\ The product $\GR_1 \times \GR_2$   
of two differentiable groupoids $\GR_1$ and $\GR_2$ is again a   
differentiable groupoid for the product structural morphisms.    
   
To obtain the Lie algebroid of the product groupoid we use   
the external product of vector bundles (not the direct sum):\    
$A(\GR_1 \times \GR_2) = A(\GR_1) \times A(\GR_2)$, a vector bundle   
over $\GR_1^{(0)}\times \GR_2^{(0)}$.   
\end{example}

We now include some examples that are suggested by the general theory of    
differentiable groupoids.

\begin{example}\ {\em Bundles of Lie groups:\ }\ In this example $\GR$ is   
a fiber bundle $p:\GR \to B$ such that each fiber   
$\GR_b := p^{-1}(b)$ has a Lie   
group structure, and the induced map    
$$   
        \GR \times_B \GR := \{ (g,g') \in \GR \times \GR,   
        p(g) = p(g') \} \ni (g,g') \mapsto g^{-1}g' \in \GR   
$$    
is a smooth map. We define $d = r = p$. The units   
of the groups $\GR_b$ then form a submanifold of $\GR$ diffeomorphic to $B$   
via $p$. We do not assume that the fibers $\GR_b$ are all    
isomorphic, but this is true in most cases of interest.

Let $\lgg$ be the restriction to the space of units of   
the vertical tangent bundle to the fibration $\GR \to B$. Then $\lgg$ is a    
vector bundle over $B$ and $A(\GR) = \lgg$. The fiber of $\lgg$ above $b$   
is then the Lie algebra of $\GR_b$, and $\lgg$ is a bundle of Lie algebras.   
In    
this   
example, the anchor map $\varrho:A(\GR) \to TB$ is the zero map.   
\end{example}

\begin{example}\ {\em Fibered products:\ }\ Let $\GR_1$ and $\GR_2$ be two    
differentiable groupoids with units $M_1$ and $M_2$. We assume that    
both $M_1$ and   
$M_2$ come equipped with submersions $p_i: M_i \to B$, $i = 1,2$, for    
some common   
manifold with corners $B$. Suppose that for each arrow    
$g \in M_i$, $p_i(d(g)) = p_i(r(g)) =:p_{i}(g)$.    
The fibered product of $\GR_1$ with $\GR_2$    
(with respect to $p_i$) is   
then    
$$   
        \GR_1 \times_B \GR_2 := \{ (g_1,g_2) \in \GR_1 \times_B \GR_2 ,    
        p_1(g_1) = p_2(g_2)\},   
$$   
with product (and structural maps, in general) induced from the product   
groupoid $\GR_1 \times \GR_2$.     
   
We get $A(\GR_1 \times_B \GR_2) = A(\GR_1) \times_B A(\GR_2)$.   
\end{example}

We are particularly interested in the above example when $\GR_1 = Y   
\times_B Y$ is the fibered pair groupoid of Example   
\ref{ex.fibered.pair} and $\GR_2$ is a bundle of Lie groups with base   
$B$. This situation seems to be fundamental in the study of   
pseudodifferential operators associated to various Lie algebras of   
vector fields.  It will be used for example in Example   
\ref{ex.boundary.f1}.  An index theorem in this framework was obtained   
in \cite{NistorIndFam}, if the fibers of $\GR_2 \to B$ are   
simply-connected solvable.

We continue with some more elaborate examples.

\begin{example}\ {\em The graph of a foliation:\ }\    
This groupoid was introduced in    
\cite{Winkelnkemper1}. Let $(M,F)$ be a foliated manifold. Thus $F \subset TM$    
is an integrable bundle. The graph of $(M,F)$ consists of equivalence classes   
$[\gamma]$ of paths $\gamma$ which are completely contained in a leaf, with   
respect to the equivalence relation    
$[\gamma] = [\gamma']$ if, and only if, $\gamma$ and    
$\gamma'$ have the same holonomy (this implies, in particular,    
that they have the    
same end-points).   
   
The Lie algebroid  is $A(\GR) = F$ in this example.   
\end{example}

\begin{example}\ {\em The fundamental groupoid:\ }\   
Let $\GR$ be the fundamental groupoid of a compact smooth manifold $M$   
(without corners) with fundamental group $\pi_1(M)=\Gamma$.  Recall   
that if we denote by $\widetilde{M}$ a universal covering of $M$ and   
let $\Gamma$ act by covering transformations on $\widetilde{M}$, then   
we have $\Gr0=\widetilde{M}/\Gamma = M$, $\GR =(\widetilde{M} \times   
\widetilde{M})/\Gamma$, and $d$ and $r$ are the two projections. Each   
fiber $\GR_x$ can be identified with $\widetilde{M}$, uniquely up to   
the action of an element in $\Gamma$.   
   
The Lie algebroid is $TM$, as in the first example.   
\end{example}

The following example generalizes the tangent groupoid of Connes; here 
we closely follow \cite[II,5]{connes}.  A similar construction was 
used in \cite{Skandalis} to define the so called ``normal groupoid,'' 
which is, anticipating a little bit, the adiabatic groupoid of a 
foliation.

\begin{example}\label{ex.adiabatic}\ {\em The adiabatic groupoid:\ }\   
The adiabatic groupoid $\GR_{\adb}$ associated to a differentiable  
groupoid $\GR$ is defined as follows. The space of units is  
$$   
        \GR_{\adb}^{(0)} := [0,\delta) \times \Gr0\,, \quad \delta > 0,   
$$   
with the product manifold structure. The set of arrows is defined as  
the disjoint union  
$$   
        \GR_{\adb} := A(\GR) \cup (0,\delta) \times \GR.  
$$   
The structural maps are defined as follows. The domain and range are:  
$$   
        d(t,g)=(t,d(g))\,\quad r(t,g)=(t,r(g))\,,\quad t>0,   
$$   
and $d(v)=r(v)=(0,x)$, if $v \in T_x\GR_x$.    
The composition is    
$   
        \mu(\gamma,\gamma')=(t,gg'),   
$   
if $t >0$, $\gamma=(t,g)$, and $\gamma'=(t,g')$,   
and    
$$   
        \mu(v,v')=v +v' \quad \text{ if } v,v' \in T_x\GR_x.   
$$    
   
The smooth structure on the set of arrows is the product structure for  
$t > 0$.  In order to define a coordinate chart at a point  
$$   
        v \in T_x\GR_x = A_x(\GR) = d^{-1}(0,x),   
$$    
choose first a coordinate system $\psi : U = U_1 \times U_2 \to \GR$,  
$U_1 \subset \RR^p$ and $U_2\subset \RR^n$ being open sets containing  
the origin, $U_2$ convex, with the following properties:\ $\psi(0,0)=x  
\in M \subset \GR$, $\psi(U) \cap M=\psi(U_1 \times \{0\})$, and there  
exists a diffeomorphism $\phi : U_1 \to \Gr0$ such that $d(\psi(s,y))  
= \phi(s)$ for all $y \in U_2$ and $s\in U_{1}$.  
   
We identify, using the differential $\operatorname{D}_2\! \psi$ of the  
map $\psi$, the vector space $\{s\} \times \RR^n$ and the tangent  
space $T_{\phi(s)}\GR_{\phi(s)}=A_{\phi(s)}(\GR)$.  We obtain then  
coordinate charts $ \psi_\varepsilon : A(\GR)\vert_{\phi(U_1)} \times  
(0,\varepsilon) \times U_1 \times \varepsilon^{-1}U_2 \to \GR, $  
$$   
        \psi_\varepsilon(0,s,y)=(0,(\operatorname{D}_2\!\psi) (s,y))  
        \in T_{\phi(s)}\GR_{\phi(s)}=A_{\phi(s)}(\GR)  
$$   
and $\psi_\varepsilon(t,s,y)=(t,\psi(s,ty)) \in (0,1) \times \GR$. For  
$\varepsilon$ small enough, the range of $\psi_\varepsilon$ will  
contain $v$.  
   
The Lie algebroid of $\GR_{\adb}$ is the adiabatic Lie algebroid  
associated to $A(\GR)$, $A(\GR_{\adb})=A(\GR)_{t}$, for all $t$, such  
that $\Gamma(A(\GR)) \cong t\Gamma(A(\GR \times [0,\delta)))$.  
\end{example}

We expect the above constructions to have applications to  
semi-classical trace formul\ae, see Uribe's overview \cite{Uribe}.  A  
variant of the above example can be used to treat adiabatic limits  
when the metric is blown up in the base. See \cite{Witten1} for some  
connections with physics.  
   
More examples are discussed in Section \ref{Sec.Examples.II}.

\section{Pseudodifferential operators on groupoids\label{Sec.POOG}}

We proceed now to define the space of pseudodifferential operators    
acting on sections of vector bundles on a differentiable groupoid.   
This construction is the same as the one in \cite{NWX}, but slightly more   
general because we consider also certain non-Hausdorff groupoids.    
General reference for pseudodifferential operators on smooth manifolds   
are, for instance,  \cite{Hormander3} or \cite{shubin}.    
We discuss operators on functions, for simplicity,   
but at the end we briefly indicate the changes necessary to handle    
operators between sections of smooth vector bundles.   
   
Our construction of pseudodifferential operators on groupoids is obtained    
considering   
certain families of pseudodifferential operators on smooth, generally   
non-compact manifolds. We begin then by recalling a few facts about    
pseudodifferential   
operators on smooth manifolds.   
   
Let $W \subset \RR^N$ be an open subset. Define the space    
${\mathcal S}^m(W \times \RR^n)$ of    
{\em symbols of order $m\in\RR$}   
on the bundle $W \times \RR^n \to W$,  as in \cite{Hormander3},    
to be the set of smooth functions $a : W \times \RR^n \to \CC$ such that   
\begin{equation}\label{eq.symbol.estimates}   
        |\partial_y^\alpha \partial_\xi^\beta a(y,\xi)| \leq C_{K,\alpha,\beta}   
        (1+|\xi|)^{m-|\beta|}   
\end{equation}   
for any compact set $K\subset W$ and any multi-indices $\alpha$ and $\beta$,    
and   
some constant $C_{K,\alpha,\beta} > 0$.   
A symbol $a \in {\mathcal S}^m(W \times \RR^n)$ is called {\em classical} if   
it has an asymptotic expansion as an infinite sum of homogeneous symbols   
$a \sim \sum_{k=0}^\infty a_{m-k}$, $a_l$ homogeneous    
of degree $l$   
for large $\norm{\xi}$, i.e.\  $a_l(y,t\xi)=t^la_l(y,\xi)$ if $\|\xi\|\geq 1$   
and $t \geq 1$.  More precisely, $\sim$ means   
$$   
        a -\sum_{k=0}^{M-1} a_{m-k}  \in     
        {\mathcal S}^{m-M}(W \times \RR^n) \quad \text{ for all } M\in\nzn\,.   
$$   
The space of classical symbols   
will be denoted by ${\mathcal S}^m_{\cl}(W \times \RR^n)$.    
Using local trivializations   
the definition of (classical)    
symbols immediately extends to arbitrary vector bundles   
$\pi:E\longrightarrow M$.   
We shall    
consider only classical symbols in this paper.   
For $a \in {\mathcal S}^m(T^*W) = {\mathcal S}^m(W \times \RR^n)$ and    
$W$  an open subset of $\RR^n$, we    
define  an operator $a(y,D_y):\CIc(W) \to \CI(W)$ by   
\begin{equation}   
        a(y,D_y)u(y)=(2\pi)^{-n}\int_{\RR^n}e^{i y\cdot \xi}   
        a(y,\xi)\hat{u}(\xi)d\xi\,,   
\end{equation}    
where $\hat{u}$ denotes the Fourier transform of $u$.   
   
Recall that if $M$ is a smooth manifold, a linear map $T:\CIc(M) \to  
\CI(M)$ is called {\em regularizing} if, and only if, it has a smooth  
distributional (or Schwartz) kernel.  Also, recall that a linear map  
$P:\CIc(M) \to \CI(M)$ is called a {\em (classical) pseudodifferential  
operator of order $m$} if, and only if, for all smooth functions  
$\phi$ supported in a (not necessarily connected) coordinate chart  
$W$, the operator $\phi P \phi$ is of the form $a(y,D_y)$ with a  
(classical) symbol $a$ of order $m$.  For a classical  
pseudodifferential operator $P$ as the one considered here, the  
collection of all classes of $a$ in ${\mathcal  
S}_{\cl}^m(T^*W)/{\mathcal S}_{\cl}^{m-1}(T^*W)$, for all coordinate  
neighborhoods $W$, patches together to define a class $\sigma_m(P)\in  
{\mathcal S}_{\cl}^m(T^*W)/{\mathcal S}_{\cl}^{m-1}(T^*W)$, which is  
called {\em the homogeneous principal symbol} of $P$; the latter space  
can, of course, canonically be identified with $S^{[m]}(T^*M)$, the  
space of all smooth functions $T^*M\setminus\{0\}\longrightarrow \CC$  
that are positively homogeneous of degree $m$. We shall sometimes  
refer to a classical pseudodifferential operator acting on a smooth  
manifold (without corners) as an {\em ordinary} classical  
pseudodifferential operator, in order to distinguish it from a  
pseudodifferential operator on a groupoid.  
   
We now begin the discussion of pseudodifferential operators on  
groupoids. A pseudodifferential operator on a differentiable groupoid  
$\GR$ will be a family $(P_x)$, $x \in M$, of classical  
pseudodifferential operators $P_x : \CIc(\GR_x) \to \CI(\GR_x)$ with  
certain additional properties that need to be explained.  
   
If $(P_x)$, $x \in M$, is a family of pseudodifferential operators acting on   
$\GR_x$, we denote by $k_{x}$ the distribution kernel of $P_x$   
We then define the support of the operator $P$ to be   
\begin{equation}\label{support}   
        \supp (P) = \overline{\bigcup_{x\in M} \supp(k_{x})}.   
\end{equation}   
The support of $P$ is contained in the closed subset   
$\{(g,g'), d(g)=d(g')\}$ of the product $\GR \times \GR$.

To define our class of pseudodifferential operators, we shall need   
various conditions on the support of our operators. We introduce the   
following terminology: a family $P=(P_x)$, $x \in M$ is    
{\em properly supported}   
if $p_i^{-1}(K) \cap \supp(P)$ is   
a compact set for any compact subset $K \subset \GR$, where   
$p_1,p_2 :\GR \times \GR \to \GR$ are the two projections.   
The family  $P = (P_x)$ is called {\em compactly supported} if its support   
$\supp(P)$ is compact. Finally, $P$ is called {\em uniformly supported} if   
its {\em reduced support} $\supp_\mu(P):=\mu_1(\supp(P))$ is    
a compact subset of $\GR$, where $\mu_1(g',g):=g'g^{-1}$.    
Clearly, a uniformly supported operator is properly supported, and a    
compactly supported operator is uniformly supported. If the family   
$\FAM$ is properly supported, then each $P_x$ is properly supported, but   
the converse is not true.   
   
Recall that the composition of two ordinary pseudodifferential operators is    
defined if one of them is properly supported. It follows that we can define    
the composition $PQ$ of two properly supported families of operators    
$P=(P_x)$ and $Q=(Q_x)$ acting on the fibers of $d : \GR \to M$   
by pointwise composition $PQ=(P_xQ_x)$,  $x \in M$. The resulting   
family $PQ$ will also be properly supported. If $P$ and $Q$ are   
uniformly supported, then $PQ$ will also be uniformly supported.     
   
The action of a family $P=(P_x)$ on functions on $\GR$ is defined pointwise    
as follows.   
For any smooth function $f \in  \CIc(\GR)$    
denote by $f_x$ its restriction $f\vert_{\GR_x}$.    
If each $f_x$ has compact support, and $\FAM$ is a family of   
ordinary pseudodifferential operators, then we define $Pf$ by    
$(Pf)_x = P_{x}(f_x).$ If $P$ is uniformly supported, then $Pf$ is   
also compactly supported. However, it is not true that $Pf \in\CIc(\GR)$   
if $f \in \CIc(\GR)$, in general, so we need some  conditions    
on the family $P$.    
We shall hence consider uniformly supported families $P = (P_x)$ because   
this guarantees that $Pf$ has compact support if $f$ does.   
    
A {\em fiber preserving diffeomorphism} will be a    
diffeomorphism    
$\psi : d(V) \times W \to V$ satisfying $d(\psi(x,w))=x$, where $W$ is some   
open subset of an Euclidean space of the appropriate dimension.   
We now discuss the differentiability condition on a family $P= (P_x)$,   
a condition which, when satisfied, implies that $Pf$ is smooth for all   
smooth $f \in \CIc(\GR)$.

\begin{definition}\label{Def.Differentiable}\   
Let $\GR$ be a differentiable groupoid with units $M$.   
A family $(P_x)$ of  pseudodifferential   
operators acting on $\CIc(\GR_x)$, $x \in M$,   
is called \underline{differentiable}   
if, and only if, for any fiber preserving diffeomorphism     
$\psi : d(V) \times W \to V$    
onto an open set $V \subseteq \GR$, and for any $\phi \in \CIc(V)$,    
we can find $a \in {\mathcal S}_{\cl}^m(d(V) \times T^*W)$   
such that $\phi P_x \phi$ corresponds to $a(x,y,D_y)$    
under the diffeomorphism  $\GR_x\cap V \simeq W$, for each $x \in d(V)$.   
\end{definition}

Thus, we require that the operators $P_x$ be    
given  in local coordinates by symbols $a_x$ that depend smoothly on   
all variables. Note that nowhere in the above definition it is necessary   
for  $\GR$ to   
be Hausdorff.  All we do need is that each of $\GR_x$ and $M = \Gr0$   
are Hausdorff.   
   
To define pseudodifferential operators on $\GR$ we shall consider   
smooth, uniformly supported families $P=(P_x)$ that satisfy also an   
invariance condition. To introduce this invariance condition, observe   
that right translations on $\GR$ define linear isomorphisms    
\begin{equation}   
\label{eq.isomorphism}   
        U_g:\CI(\GR_{d(g)}) \to \CI(\GR_{r(g)}):   
        (U_gf)(g')=f(g'g) \,.   
\end{equation}   
A family of operators $P=(P_x)$ is then called {\em invariant}   
if $P_{r(g)}U_g = U_g P_{d(g)}$, for all $g \in \GR$.   
We are now ready  to define pseudodifferential operators on $\GR$.

\begin{definition}\label{Main.definition}   
Let $\GR$ be a differentiable groupoid with units $M$, and   
let $P=(P_x)$ be a family  $P_x : \CIc(\GR_x) \to \CI(\GR_x)$ of    
(order $m\in\RR$, ordinary) classical pseudodifferential operators.   
Then $P$ is an (order $m$) \underline{pseudodifferential operator on    
$\GR$} if, and only if, it is   
   
(i)\  uniformly supported,   
   
(ii)\ differentiable, and    
   
(iii)\ invariant.   
   
\noindent   
We denote the space of order $m$ pseudodifferential operators on $\GR$   
by $\tPS{m}$.   
\end{definition}

We also denote $\tPS{\infty} := \cup_{m \in \RR}\tPS{m}$ and    
$\tPS{-\infty}=\cap_{m \in \RR}\tPS{m}$.    
Let us now give an alternative description of $\tPS{m}$ that highlights   
again the conormal nature of  kernels of pseudodifferential    
operators. For $P\in\tPS{\infty}$, we call   
$$\kappa_{P}(g):=k_{d(g)}(g,d(g))\,, g\in\GR$$ the {\em reduced} or  
{\em convolution kernel} of $P$.  Due to the right-invariance of $P$  
we expect the reduced kernel to carry all information of the family  
$P$. For the definition of conormal distributions we refer the reader  
to \cite{fio, Hormander3} in the smooth case.

\begin{proposition}   
\label{redker}\   
The map $P\longmapsto \kappa_{P}$ induces an isomorphism   
$$\tPS{m}\stackrel{\cong}{\longrightarrow}I_{c}^{m}(\GR,M;d^{*}\mathcal{D})$$   
where $I_{c}^{m}(\GR,M;d^{*}\mathcal{D})$ denotes the space of all   
compactly   
supported,   
$d^{*}\mathcal{D}$-valued distributions on $\GR$ conormal   
to $M$. In particular,   
$P\longmapsto\kappa_{P}$ identifies $\tPS{-\infty}$ with the    
convolution algebra $\CIc(\GR,d^{*}\mathcal{D})$. Moreover, we have   
$\supp(\kappa_{P})=\supp_{\mu}(P)$.   
\end{proposition}

Define $\CIc(\GR)$ as in Section    
\ref{Sec.prelim}, then $\tPS{m}(\CIc(\GR)) \subset \CIc(\GR)$.   
We obtain in this way a representation    
$\pi$ of $\tPS{\infty}$ on $\CIc(M)$, uniquely determined by   
\begin{equation}\label{eq.vector.rep}    
        (\pi(P)f) \circ r = P (f \circ r), \quad P = (P_x) \in \tPS{m}.   
\end{equation}   
We call this representation acting on any space of functions on $M$ on which   
it makes sense ($\CIc(M)$, $\CI(M)$, $L^2(M)$, or Sobolev   
spaces) the {\em vector representation} of $\tPS{\infty}$.   
   
We now discuss the extension of the principal symbol map to $\tPS{m}$.   
Denote by  $\pi:A^*(\GR)\rightarrow M$ the projection.    
If $P=(P_x)\in \tPS{m}$ is an order $m$    
pseudodifferential  operator on $\GR$, then the principal   
symbol $\sigma_m(P)$ of $P$ will be an order $m$ homogeneous function    
on $A^*(\GR)\smallsetminus 0$ (it is defined only outside the zero   
section) such that:    
\begin{equation}\label{princ.symb}   
        \sigma_m(P)(\xi)=\sigma_m(P_x)(\xi) \in \CC    
        \quad \text{ if } \xi \in A^*_x(\GR)=T^*_x\GR_x.   
\end{equation}   
Denote by ${\mathcal S}_{c}^{m} (A^*_x(\GR)) \subset {\mathcal S}_{\cl}^{m}    
(A^*_x(\GR))$ the subspace of classical symbol whose support has compact   
projection onto the space of units $M$.   
The above definition determines a linear map    
\begin{equation*}   
        \sigma_m:\tPS{m} \to {\mathcal S}_{c}^m(A^*(\GR))/   
        {\mathcal S}_{c}^{m-1}(A^*(\GR))   
\cong S^{[m]}_{c}(A^{*}(\GR))\,,   
\end{equation*}   
where $S^{[m]}_{c}(A^{*}(\GR))$ denotes the space of   
all smooth functions    
$A^{*}(\GR)\setminus\{0\}\longrightarrow\CC$   
that are positively homogeneous of degree $m$, and whose   
support has compact projection   
onto the space of units $M$.   
The map $\sigma_{m}$ is said to be the    
{\em principal symbol}. A   
pseudodifferential operator $P\in\tPS{m}$ is called   
{\em elliptic} provided its principal   
symbol $\sigma_{m}(P)\in  S^{[m]}_{c}(A^{*}(\GR))$   
does not vanish on $A^{*}(\GR)\setminus\{0\}$; note,   
with this definition   
elliptic operators only exist if   
the space of units is compact.

The following result extends several of the well-known properties of the    
calculus   
of ordinary pseudodifferential operators on smooth manifolds to   
$\tPS{\infty}$. Denote by $\{\; , \; \}$ the canonical Poisson   
bracket on $A^*(\GR)$.

\begin{theorem}\label{Theorem.Algebra}    
Let $\GR$ be a differentiable groupoid. Then $\tPS{\infty}$ is an    
algebra    
with the following properties:    
   
(i) The principal symbol map    
$$   
        \sigma_m:\tPS{m} \to {\mathcal S}_{c}^m(A^*(\GR))/   
        {\mathcal S}_{c}^{m-1}(A^*(\GR))   
$$    
is surjective with kernel $\tPS{m-1}$.   
   
(ii) If $P \in \tPS{m}$ and $Q \in \tPS{m'}$, then $PQ \in \tPS{m+m'}$  
and satisfies  
$\sigma_{m+m'}(PQ)=\sigma_m(P)\sigma_{m'}(Q)$. Consequently, $[P,Q] \in  
\tPS{m + m' -1}$. Its principal symbol is given by $\sigma_{m + m'  
-1}([P,Q]) = \frac{1}{i}\{ \sigma_m(P), \sigma_{m'}(Q) \}$.  
\end{theorem}

Properly supported invariant differentiable families of  
pseudodifferential operators also form a filtered algebra, denoted  
$\Psi_{\prop}^{\infty}(\GR)$.  While it is clear that, in order for  
our class of pseudodifferential operators to form an algebra, we need  
some condition on the support of their distributional kernels, exactly  
what support condition to impose is a matter of choice. We prefer the  
uniform support condition because it leads to a better control at  
infinity of the family of operators $P=(P_x)$ and allows us to  
identify the regularizing ideal (\ie the ideal of order $-\infty$  
operators) with the groupoid convolution algebra of $\GR$. The choice  
of uniform support will also ensure that $\tPS{m}$ behaves  
functorially with respect to open embeddings. The compact support  
condition enjoys the same properties but is usually too  
restrictive. The issue of support will be discussed again in examples.  
   
We now discuss the restriction of families in $\tPS{m}$ to invariant  
subsets of $M$, or, more precisely, the restriction to $\GR_Y$, the  
``reduction'' of $\GR$ to the invariant subset $Y$ of $M$. The  
resulting restriction morphisms $\inn_Y : \tPS{\infty} \to  
\Psi^\infty(\GR_Y)$ is the analog in our setting of the indicial  
morphisms considered in \cite{mepi92}.  
    
We continue to denote by $\GR$ a differentiable groupoid with units  
$M$.  Let $Y \subset M$ and let $\GR_Y := d^{-1}(Y) \cap r^{-1}(Y)$.  
Then $\GR_Y$ is a groupoid with units $Y$, called the {\em reduction  
of $\GR$ to $Y$}. An {\it invariant} subset $Y \subset M$ is a subset  
such that $d(g) \in Y$ implies $r(g) \in Y$. For an invariant subset  
$Y \subset M$, the reduction of $\GR$ to $Y$ satisfies  
$$   
        \GR_Y = d^{-1}(Y) = r^{-1}(Y)   
$$    
and is a differentiable groupoid, if $Y$ is a closed submanifold (with  
corners) of $M$.  
   
If $P = (P_x)$, $x \in M$, is a pseudodifferential operator on   
$\GR$, and $Y \subset M$ is an closed, invariant submanifold with corners,    
we can restrict $P$ to $d^{-1}(Y)$ and obtain   
$$   
        \inn_Y (P) := (P_x)_{x\in Y},   
$$   
which is a family of operators acting on the fibers of   
$d:\GR_Y=d^{-1}(Y) \to Y$ and satisfying all the conditions necessary   
to define an element of $\Psi^{\infty}(\GR\vert_Y)$.  This leads to a map   
\begin{equation}\label{eq.indicial}    
        \inn_Y =\inn_{Y,M}:\tPS{\infty} \to \Psi^{\infty}(\GR\vert_Y),   
\end{equation}   
which is easily seen to be an algebra morphism.    
   
Let us indicate now what changes need to be made when we consider operators   
acting on sections of a vector bundles. Because operators acting    
between sections    
of two {\em different} vector bundles $E_1$ and $E_2$ can be recovered from    
operators acting on $E=E_1 \oplus E_2$, we may assume that $E_1 = E_2 = E$ as   
vector bundles on $M = \Gr0$.   
   
Denote by $r^*(E)$ the pull-back of $E$ to $\GR = \Gr1$.    
Then the isomorphisms of Equation \eqref{eq.isomorphism}    
will have to be replaced by   
 $$       U_g:\CI(\GR_{d(g)},r^*(E)) \to \CI(\GR_{r(g)},r^*(E)):   
        (U_gf)(g')=f(g'g) \in (r^*E)_{g'},   
$$   
which  makes sense because of $(r^*E)_{g'}=(r^*E)_{g'g}=E_{r(g')}$.     
Then,  to define $\PS{m}$ we consider families   
$P = (P_x)$ of order $m$ pseudodifferential operators $P_x$, $x \in M$,   
acting on the spaces $\CIc(\GR_x,r^*(E))$. We require these families to   
be uniformly supported, differentiable, and invariant, as in the   
case $E= \CC$.    
   
The principal symbol $\sigma_m(P)$ of a classical pseudodifferential  
operator $P$ belongs then to ${\mathcal  
S}_{c}^m(T^*W;\End(E))/{\mathcal S}_{c}^{m-1}(T^*W;\End(E))$.  
Finally, the restriction (or indicial) morphism is a morphism.  
$$   
        \inn_Y :\Psi^{\infty}(\GR;E) \to\Psi^{\infty}   
        (\GR\vert_Y ; E\vert_Y).   
$$     
All the other changes needed to treat the case of    
non-trivial vector bundles $E$    
are similar.   
   
There is one particular case of a bundle $E$ that deserves   
special attention.    
Let  $E := \VD^{1/2}$ be the square root of the density bundle    
$\VD = | \Lambda^n A(\GR)|,$   
as before.  If $P \in\Psi^{m}(\GR; \VD^{1/2})$    
consists of the family $(P_x, x\in M)$, then   
each $P_x$ acts on $V_x = C_c^{\infty}(\GR_x;r^*(\VD^{1/2})).$   
Since $r^*(\VD^{1/2}) = \Omega_{\GR_x}^{1/2}$ is the bundle of   
half densities on $\GR_x$, we can define a natural hermitian    
inner product $(\; ,\; )$ on $V_x$.    
Let $P = (P_x) \in \Psi^m(\GR;\VD^{1/2})$,    
 and denote by $P_x^*$ the formal adjoint of $P_x$, that   
is, the unique pseudodifferential operator on $V_x$    
such that $(P_x f, g) = (f, P_x^*g)$,    
for all $f,g \in V_x$. It is not hard to see that    
$(P_x^*) \in \Psi^m(\GR;\VD^{1/2})$.     
Moreover, $\sigma_m(P^*) = \overline{\sigma_m(P)}$.   
   
If $E$ is the complexification of a real bundle    
$E_0$:\ $E \simeq E_0 \otimes \CC$,   
then the complex conjugation operator $J \in \End_\RR(E)$ defines a real structure    
on $\PS{*}$, that is, a conjugate linear involution on $\PS{*}$.   
In this case, $\PS{*}$ is the complexification of the set of fixed points of   
this involution.

\section{Bounded representations\label{Sec.BR}}

For  a smooth, compact manifold $M$ (without corners),   the algebra    
$\Psi^0(M)$ of order zero pseudodifferential operators on $M$ acts by   
bounded operators on $L^2(M,d\mu)$, the Hilbert space $L^2(M,d\mu)$ being    
defined with respect to the (essentially unique) measure $\mu$ corresponding    
to a nowhere vanishing density on $M$. Moreover, this is basically the only   
interesting $*$-representation of $\Psi^0(M)$ by bounded operators on an    
infinite-dimensional Hilbert space of functions.   
   
For a differentiable groupoid $\GR$ with units $M$, a manifold with  
corners, it is still true that we can find a measure $\mu$ such that  
$\tPS0$ acts by bounded operators on $L^2(M,d\mu)$.  However, in this  
case there may exist natural measures $d\mu$ that are singular with  
respect to the measure defined by a nowhere vanishing  
density. Moreover, there may exist several non-equivalent such  
measures, and these representations may not exhaust all equivalence  
classes of non-trivial, irreducible, bounded representations of  
$\tPS{0}$.  
    
The purpose of this section is to introduce the class of  
representations we are interested in, and to study some of their  
properties. A consequence of our results is that in order to construct  
and classify bounded representations of $\tPS0$, it is essentially  
enough to do this for $\tPS{-\infty}$.  
   
We are interested in representations of $\tPS{m}$, $m \in  
\{0,\pm\infty\}$.  We fix a trivialization of $\VD$, so that we get an  
isomorphism $\tPS{m} \cong \Psi^m(\GR;\VD^{1/2})$ and hence we have an  
involution $*$ on $\tPS{m}$.  Let ${\mathcal H}_0$ be a dense subspace  
of a Hilbert space $\mathcal{H}$ with inner product $(\; , \;  
)$. Recall that a $*$-morphism $\alpha : A \to \End(\mathcal H_0)$  
from a $*$-algebra $A$ is a morphism such that $(\alpha(P^*)\xi, \eta)  
= (\xi, \alpha(P)\eta))$, for all $P \in A$ and all $\xi,\eta \in  
\mathcal H_0$.

\begin{definition}\label{def.bounded.rep}\   
Let $\mathcal{H}_{0}$ be a dense subspace of a  
Hilbert space $\mathcal{H}$, and $m = 0$ or $m =    
\pm \infty$ be fixed.  
A \underline{bounded representation} of $\tPS{m}$ on   
${\mathcal H}_0$ is a $*$-morphism   
$\varrho : \tPS{m} \to \End( {\mathcal H}_0 )$ such that   
$\varrho(P)$ extends to a bounded   
operator on $\mathcal{H}$  
for all $P \in \tPS{0}$ (for all $P \in \tPS{-\infty}$,    
if $m = -\infty$).   
\end{definition}

Note that if $\varrho$ is as above and $P$ is an operator of positive  
order, then $\varrho(P)$ does not have to be bounded.  Since  
$\mathcal{H}_{0}$ is dense in $\mathcal{H}$, each operator  
$\varrho(P)$ can be regarded as a densely defined operator.  The  
definition implies that $\varrho(P^*) \subset \varrho(P)^*$, so the  
adjoint of $\varrho(P)$ is densely defined, and hence $\varrho(P)$ is  
a closable operator.  We shall usually make no distinction between  
$\varrho(P)$ and its closure.  
  
We call a bounded representation  
$\varrho:\tPS{-\infty}\rightarrow\End({\mathcal H})$ {\em  
non-degenerate} provided $\varrho(\tPS{-\infty})\mathcal{H}$ is dense  
in $\mathcal{H}$.  
   
The following theorem establishes, among other things, a bijective  
correspondence between non-degenerate bounded representations of  
$\tPS{-\infty}$ on a Hilbert space $\mathcal{H}$, and bounded  
representations $\varrho : \tPS{\infty} \rightarrow\End(\mathcal H_0)$  
such that the space $\varrho(\tPS{-\infty})\mathcal H_0$ is dense in  
$\mathcal{H}$. We shall need the following slight extension of a  
result in \cite{LMN}.

\begin{theorem} \label{Theorem.EXT}\   
Let ${\mathcal H}$ be a Hilbert space and let $\varrho : \tPS{-\infty} \to   
\End({\mathcal H})$ be a bounded representation. Then, to each  
$P \in \tPS{s}$, $s \in \RR$, we can associate an unbounded operator 
$\varrho(P)$, with domain ${\mathcal H_0}:= \varrho(\tPS{-\infty}) 
{\mathcal H}$, such that $\varrho(P)\varrho(R) = 
\varrho(PR)$ and $\varrho(R) \varrho(P) = \varrho(RP)$, for any $R \in  
\tPS{-\infty}$.  
  
We obtain in this way an extension of $\varrho$   
to a bounded representation of $\tPS{0}$ on ${\mathcal H}$ and to a  
bounded representation of $\tPS{\infty}$ on ${\mathcal H_0}:=  
\varrho(\tPS{-\infty}){\mathcal H}$.   
\end{theorem}

\begin{proof}\ We may assume that  
${\mathcal H}_0$ is dense in ${\mathcal H}$. Fix $P \in \tPS{0}$ and 
we let $$ \varrho(P)\xi=\varrho(PQ)\eta, $$ if $\xi = \varrho(Q)\eta$, 
for some $Q \in \Psi^{-\infty}(\GR)$ and $\eta 
\in {\mathcal H}$. We need to show that $\varrho(P)$ is well-defined and  
bounded. Thus, we need to prove that $\sum_{k=1}^N \varrho(PQ_k)\xi_k=0$,  
if $P \in \tPS{0}$ and $\sum_{k=1}^N \varrho(Q_k)\xi_k=0$, for some $Q_k  
\in \tPS{-\infty}$ and $\xi_k \in {\mathcal H}$.  
  
We will show that, for each $P \in \tPS0$, there exists a constant  
$k_P>0$ such that  
\begin{equation}  
        \| \sum_{k=1}^N\varrho(PQ_k)\xi_k \| \le k_P   
        \| \sum_{k=1}^N\varrho(Q_k)\xi_k \|.  
\end{equation}  
Let $C \ge | \sigma_0(P)| + 1$, and let  
\begin{equation}\label{eq.b}  
        b = (C^2 - |\sigma_0(P)|^2)^{1/2}.  
\end{equation}  
Then $b-C$ is in $\CIc(S^*(\GR))$, and it follows from Theorem  
\ref{Theorem.Algebra} that we can find $Q_0 \in\tPS0$ such that  
$\sigma_0(Q_0) = b-C$. Let $Q = Q_0 + C$.  Using again Theorem  
\ref{Theorem.Algebra}, we obtain, for  
$$  
        R = C^2 - P^*P - Q^*Q \in \tPS0,  
$$   
that  
$$  
        \sigma_0(R) = \sigma_0(C^2 - P^* P - Q^* Q)= 0,  
$$  
and hence $R \in\tPS{-1}$. A standard argument using the  
asymptotic completeness of the algebra of pseudodifferential  
operators shows that we can assume that $Q$ has  
order $-\infty$. Let then  
\begin{equation}\label{eq.xieta}  
        \xi = \sum_{k=1}^N \varrho(Q_k)\xi_k, \;\,  
        \eta = \sum_{k=1}^N \varrho(PQ_k)\xi_k, \quad \mbox{\rm and } \quad  
        \zeta= \sum_{k=1}^N \varrho(QQ_k)\xi_k,  
\end{equation}  
which gives,  
\begin{multline}\label{eq.middle}  
        \|\eta\|^2=(\eta,\eta)=\sum_{j,k=1}^N  
        (\varrho(Q_k^*P^*PQ_j)\xi_j,\xi_k)\\ = \sum_{j,k=1}^N \big(C^2  
        (\varrho(Q_kQ_j)\xi_j,\xi_k) - (\varrho(Q_k^*Q^*QQ_j)\xi_j,\xi_k) -  
        (\varrho(Q_k^*RQ_j)\xi_j,\xi_k)\big ) \\ = C^2 \|\xi\|^2 -  
        \|\zeta\|^2 - (\varrho(R)\xi,\xi) \leq (C^2 + \|\varrho(R)\|)  
        \|\xi\|^2.  
\end{multline}  
  
The desired representation of $\tPS0$ on ${\mathcal H}$ is obtained by  
extending $\varrho(P)$ by continuity to ${\mathcal H}$.  
  
To extend $\varrho$ further to $\tPS{s}$, we proceed similarly: we want  
$$  
        \varrho(P)\varrho(Q)\xi=\varrho(PQ)\xi,  
$$  
for $P \in \tPS{\infty}$ and $Q \in \tPS{-\infty}$. Let $\xi$ and  
$\eta$ be as in Equation \eqref{eq.xieta}. We need to prove that $\eta  
= 0$ if $\xi = 0$.  Now, because ${\mathcal H}_0$ is dense in  
${\mathcal H}$, we can find $T_{j}$ in $A_\varrho$ the norm closure of  
$\varrho(\tPS{-\infty})$ and $\eta_{j} \in {\mathcal H}$ such that $\eta  
= \sum_{j=1}^{N}T_{j}\eta_{j}$.  Choose an approximate unit $u_\alpha$  
of the $C^*$-algebra $A_\varrho$, then $u_\alpha T_{j} \to T_{j}$ (in the  
sense of generalized sequences). We can replace then the generalized  
sequence (net) $u_\alpha$ by a subsequence, call it $u_m$ such that  
$u_m T_{j} \to T_{j}$, as $m \to \infty$. By density, we may assume  
$u_m = \varrho(R_m)$, for some $R_m \in \tPS{-\infty}$. Consequently,  
$\varrho(R_m )\eta \to \eta$, as $m \to \infty$. Then  
$$  
        \eta = \lim \sum_{k = 1}^N \varrho(R_m) \varrho(PQ_k) \xi_k     
        = \lim \sum_{k = 1}^N \varrho(R_m P)\varrho(Q_k) \xi_k = 0,  
$$  
because $R_mP \in \tPS{-\infty}$.   
\end{proof}  
  
\vspace*{1mm} {\em Remark.}\ We also obtain, using the above notation,  
that any extension of $\varrho$ to a representation of $\tPS0$ is  
bounded.  This extension is uniquely determined if ${\mathcal H_0}$ is  
dense in ${\mathcal H}$. \vspace*{1mm}

Assume that $M$ is connected, so that the manifolds $\GR_x$ have the  
same dimension $n$.  We now proceed to define a Banach norm on  
$\tPS{-n-1}$.  This norm depends on the choice of a trivialization of  
the bundle of densities $\VD$, which then gives rise to a right  
invariant system of measures $\mu_x$. Then, if $P \in \tPS{-n-1}$, we  
use the chosen trivialization of $\VD$ to identify the reduced kernel  
$\kappa_P$, which is a priori a distribution, with a compactly  
supported, {\em continuous function} on $\GR$, still denoted by  
$\kappa_P$.  We then define  
\begin{equation}\label{eq.norm.one}   
        \|P\|_1 =\sup_{x \in M}\left\{\int_{\GR_x}   
        |\kappa_P(g^{-1})|d\mu_x(g), \int_{\GR_x}   
        |\kappa_P(g)|d\mu_x(g)\right\}.   
\end{equation}

Some of the most interesting representations of $\tPS{\infty}$ are the  
regular representations $\pi_x$, $x \in M$. These are bounded  
representations defined as follows; let $x \in M$, then the {\em  
regular representation} $\pi_x$ associated to $x$ is the natural  
representation of $\tPS{\infty}$ on  
$C_c^{\infty}(\GR_x;r^*(\VD^{1/2}))$, that is  
$\pi_x(P)=P_x$. Moreover, $\|\pi_x(P)\| \le \|P\|_1$, if $P \in  
\tPS{-n-1}$.  
    
Define now the {\em reduced norm} of $P$ by  
$$   
        \|P\|_r = \sup_{x\in M}  \| \pi_x(P)\|=\sup_{x\in M} \|P_x\|\, ,   
$$   
and the {\em full norm} of $P$ by  
$$   
        \|P\| = \sup_\varrho \| \varrho(P)\|,   
$$   
where $\varrho$ ranges through all bounded representations $\varrho$  
of $\tPS0$ satisfying  
$$   
        \|\varrho(P)\| \le \|P\|_1 \quad \text{ for all } P \in  
        \tPS{-\infty} \, .  
$$   
The above comments imply, in particular, that $\| P \|_r \le  
\|P\|$. If we have equality, we shall call $\GR$ {\em amenable},  
following the standard usage.  
   
Denote by $\alg\GR$ [respectively, by $\ralg\GR$] the closure of  
$\tPS0$ in the norm $\|\;\|$ [respectively, in the norm $\|\;\|_r$].  
Also, denote by $\ideal\GR$ [respectively, by $\rideal\GR$] the  
closure of $\tPS{-\infty}$ in the norm $\|\;\|$ [respectively, in the  
norm $\|\;\|_r$]. The principal symbol $\sigma_0$ extends by  
continuity to $\alg\GR$ and $\ralg\GR$.  
   
Let $S^*(\GR):=(A^{*}(\GR)\setminus\{0\})/\RR_{+}$ be the space of  
rays in $A^*(\GR)$.  (By choosing a metric on $A(\GR)$, we may  
identify $S^*(\GR)$ with the subset of vectors of length one in  
$A^*(\GR)$.)  Then we have the following two exact sequences of  
$C^*$-algebras:  
\begin{eqnarray*}   
     && 0 \to \ideal\GR \to \alg\GR \to \CO(S^*(\GR)) \to 0\quad  
      \text{ and }\\ && 0 \to \rideal\GR \to \ralg\GR \to  
      \CO(S^*(\GR)) \to 0\,.  
\end{eqnarray*}   
In particular, $\tPS{-\infty}$ is dense in $\tPS{-1}$.  
   
Let $Y \subset M$ be a closed, invariant submanifold with  
corners. Then we also have exact sequences  
\begin{eqnarray}   
\label{eq.seq2}   
&&        0 \to \ideal{\GR_{M \setminus Y}} \to \ideal{\GR} \to    
        \ideal{\GR_Y} \to 0\quad \text{ and }\\   
\label{eq.seq2'}   
&&        0 \to \alg{\GR_{M \setminus Y}} \to \alg{\GR} \to    
        \alg{\GR_Y} \to 0.   
\end{eqnarray}   
(We will not use that, but it is interesting to mentioned that it is  
known that there are no such exact sequence for reduced  
$C^*$-algebras, in general.)  
   
All the morphisms of the above four exact sequences are compatible  
with the complex conjugation on these algebras.

\begin{definition}\label{def.filtration}\ An \underline{invariant filtration}   
$Y_0 \subset Y_1 \subset \dots \subset Y_n = M$ is an increasing  
sequence of closed, invariant subsets of $M$ with the property that  
the closure $\overline{S}$ of each connected component $S$ of $Y_k  
\smallsetminus Y_{k-1}$ is a submanifold with corners of $M$ and that  
$\overline{S} \smallsetminus S$ is the union of the hyperfaces of  
$\overline{S}$ (that is, $S = \overline{S} \smallsetminus \pa S$).  
\end{definition}

The exact sequences defined before then give the following result:

\begin{theorem}\label{Theorem.CS}\     
Let $\GR$ be a differentiable groupoid with space of   
units $M$,  and let $\emptyset=:Y_{-1}\subset   
Y_0 \subset Y_1 \subset \dots \subset Y_n = M$ be an   
invariant filtration of $M$. Define    
${\mathfrak I}_k := \ideal{\GR_{M \smallsetminus    
Y_{k-1}}}$. Then we get a composition series of   
closed ideals   
$$   
        (0) \subset {\mathfrak I}_n \subset {\mathfrak I}_{n-1} \subset    
        \ldots \subset {\mathfrak I}_0 = \ideal\GR \subset \alg\GR\,,   
$$   
whose subquotients are determined by    
\begin{eqnarray*}   
        \sigma_{0}: \alg\GR /{\mathfrak I}_0   
        &\overset{\cong}{\lra}&\CO(S^*\GR)\,, \mbox{ and }\\   
        {\mathfrak I}_{k}/{\mathfrak I}_{k+1} &\cong&   
        \ideal{\GR_{Y_{k} \smallsetminus Y_{k-1}}}  \quad    
\mbox{ for }   0\leq k\leq n \,.   
\end{eqnarray*}   
\end{theorem}

A completely analogous result holds for the norm closure of the algebra   
$\PS{0}$, for any Hermitian vector bundle $E$. In fact, we can find an   
orthogonal projection $p_E \in M_N(\CI(M))$, for some large $N$, such   
that $E \cong p_E(M \times \CC^N)$, and hence    
$\PS{0} \cong p_E M_N(\tPS{0}) p_E$.   
   
The definition of an invariant filtration given in this paper is slightly    
more general   
than the one in \cite{LMN}, however, these definitions    
are equivalent if each $\GR_x$   
is connected. Thus, in order to avoid some unnecessary    
technical complications,    
{\em we shall assume from now on that all the fibers    
$\GR_x := d^{-1}(x)$ of $d$    
are connected}. (Recall that a   
groupoid with this property is called {\em $d$-connected}.)

We observe then, that if $(Y_k)_k$ is an invariant filtration, then each    
connected component of $Y_k \smallsetminus Y_{k-1}$ is an invariant    
subset of $M$ and    
$$   
        \ideal{\GR_{Y_k \smallsetminus Y_{k-1}}} \cong   
        \bigoplus_S \ideal{\GR_S}   
$$   
where $S$ ranges through the set of open components of    
$Y_{k} \smallsetminus Y_{k-1}$.   
Moreover, a completely similar direct sum decomposition exists for   
$\alg{\GR_{Y_k \smallsetminus Y_{k-1}}}$.   
   
A first consequence of the above theorem    
is that if each $\GR_S$ is an   
amenable groupoid (that is, $\ideal{\GR_S} \cong \rideal{\GR_S}$), then $\GR$   
is also amenable. This can be seen as follows. Using an argument based on   
induction, it is enough to prove that   
if a groupoid $\GR$ has an open invariant subset $ \mathcal O $ such that   
both $\GR_{\mathcal O}$ and $\GR_{M \smallsetminus \mathcal O}$ are   
amenable, then $\GR$ is amenable. To prove this,   
let $I$ be the kernel of the natural map  
$\ideal\GR\rightarrow\rideal\GR$ which  
is onto because its range is closed and contains the  
dense subspace $\tPS{-\infty}$. Since  
$\GR_{M\setminus \mathcal{O}}$ is amenable,  
$I$ is in the kernel of the restriction homomorphism  
$\ideal\GR\rightarrow\ideal{\GR_{M\setminus\mathcal{O}}}$,  
\ie $I$ is a subset of $\ideal{\GR_{\mathcal{O}}}$.  
But the maps $\ideal{\GR_{\mathcal O}} \to \rideal{\GR_{\mathcal O}}$  
and $\rideal{\GR_{\mathcal O}} \to \rideal{\GR}$ are both injective.  
Hence $I = 0$.

The above theorem leads to a characterization of compactness and  
Fredholmness for operators in $\tPS0$. This characterization is  
similar, and it actually contains as a particular case, the  
characterization of Fredholm operators in the ``$b$-calculus'' or one  
of its variants on manifolds with corners, see  
\cite{mepi92}. Characterizations of compact and Fredholm operators  
on manifolds with singularities were, for instance, also obtained in  
\cite{defr,man95,mame99,Melrose-Nistor1,plam86, plamsen94,braun,  
blau}.  
   
The significance of Theorem \ref{Theorem.CS} is that often   
in practice we can find nice invariant stratifications $M = \bigcup S$   
for which the subquotients $\ideal{\GR_S}$ have a relatively simpler     
structure than that of $\ideal{\GR}$ itself. An example is the    
$b$-calculus and its generalizations, the $c_n$-calculi, which are   
discussed in Section \ref{Sec.Examples.III}.   
   
In the following, we shall denote by $\otimes_{min}$ the {\em minimal}  
tensor product of $C^*$--algebras, defined using the tensor product of  
Hilbert spaces, see \cite{Sakai}.  More precisely, assume that $A_i$,  
$i = 1,2$, are $C^*$-algebras, which we may assume to be closed  
subalgebras of the algebras of bounded operators on some Hilbert  
spaces ${\mathcal H}_i$. Then the algebraic tensor product $A_1  
\otimes A_2$ acts on (the Hilbert space completion of) ${\mathcal H}_1  
\otimes {\mathcal H}_2$, and $A_1 \otimes_{min} A_2$ is defined to be  
the completion of $A_1 \otimes A_2$ with respect to the induced  
norm. The following result is sometimes useful.

\begin{proposition}\label{prop.tens}\ If $\GR_i$, $i = 0,1$, are two   
differential groupoids, then    
$$   
        \rideal{\GR_0 \times \GR_1} \simeq   
        \rideal{\GR_0} \otimes_{min} \rideal{\GR_1}.   
$$   
\end{proposition}

\section{Examples II: Pseudodifferential Operators   
\label{Sec.Examples.II}}

The examples of differentiable groupoids of Section \ref{Sec.Examples.I}   
also lead to interesting algebras of pseudodifferential operators. Many   
well-known algebras of pseudodifferential operators are in fact      
(isomorphic to)   
algebras of the form $\tPS{\infty}$. This leads to new insight    
into the structure   
of these algebras. In additions to these well-known algebras, we    
also obtain algebras    
that are difficult to describe directly, without using groupoids.    
Moreover, some   
of these algebras were not considered before the groupoids were introduced    
into the picture, nevertheless,    
these algebras are expected to play an important role in the analysis on   
certain    
non-compact manifolds.   
   
Our examples will follow in the beginning the same order as the examples    
considered in Section \ref{Sec.Examples.I}.

\begin{example}\ If $\GR = M$ is a manifold (possibly with corners),    
then we have $\tPS{\infty} \simeq    
\CIc(M)$ and $\Psi^\infty_{prop}(\GR) = \CI(M)$.   
\end{example}

Denote by $\Psi^m_{\prop}(M)$ the space of {\em properly supported}    
pseudodifferential operators on a smooth manifold $M$.

\begin{example}\ If $\GR=G$ is a Lie group,   
then $\tPS{m} \simeq \Psi^m_{\prop}(G)^G$, the algebra of properly   
supported pseudodifferential operators on $G$, invariant   
with respect to right translations. In this example, every invariant   
properly supported operator is also uniformly supported.     
\end{example}

The following example shows that the algebras of pseudodifferential operators    
(with appropriate support conditions for the Schwartz   
kernels) on a smooth manifolds without corners can be recovered as algebras    
of pseudodifferential operators on    
the pair groupoid. Let $\Psi_{\comp}^m(M)$ be the space of pseudodifferential   
operators on $M$ with compactly supported Schwartz kernels.

\begin{example}\ Suppose now that $\GR = M \times M$, with $M$ a smooth   
manifold without corners, is the pair groupoid. Then $\tPS{m} \cong     
\Psi_{comp}^m(M)$ and $\Psi_{prop}^m(\GR) \cong \Psi_{prop}^m(M)$.   
   
Moreover, the vector representation of $\tPS{\infty}$ on $\CIc(M)$   
recovers the usual action of pseudodifferential operators on functions on  
$M$.   
(Recall from   
\eqref{eq.vector.rep}, that    
the vector representation $\pi$ of $\tPS{\infty}$ is given by    
$(\pi(P)f) \circ r   
= P (f \circ r)$.)   
\end{example}

The fibered pair recovers families of operators.

\begin{example}\ If $\GR = M \times_B M$ is the fibered pair groupoid,   
for some submersion $M \to B$, then $\tPS{m}$ consists of families of   
pseudodifferential operators along the fibers of $M \to B$ such that their   
reduced kernels are compactly supported (as distributions on $\GR$).

The vector representation $\pi$ of $\tPS{\infty}$ on $\CIc(M)$ is just the   
usual action of families of pseudodifferential operators on functions,   
the action being defined fiberwise.   
\end{example}

The following three examples of algebras were probably considered only in   
the framework of groupoid algebras, although particular cases have been    
investigated before.

\begin{example}\    
\label{prgp}   
For a product groupoid $\GR = \GR_1 \times \GR_2$ there   
is no obvious description of $\Psi^\infty(\GR_1 \times \GR_2)$ in terms    
of $\Psi^\infty(\GR_1)$ and $\Psi^\infty(\GR_2)$, in general. However,    
when $\GR_1 = M_1$ is a manifold with corners (so $\GR_1$ has no   
non-trivial arrows), then $\Psi^m(\GR_1 \times \GR_2)$ consists of   
families of operators in $\Psi^m(\GR_2)$ parameterized by $M_1$.   
For smoothing operators the situation is simpler:\   
$\Psi^{-\infty}(\GR_1 \times \GR_2)$ contains naturally the tensor   
product $\Psi^{-\infty}(\GR_1) \otimes \Psi^{-\infty}(\GR_2)$ as a dense   
subset.    
   
An interesting particular case is when $\GR_1 = M \times M$, the pair  
groupoid, and $\GR_2 = \RR^q$ (that is, the groupoid associated to the  
Lie group $\RR^q$), then $\Psi^m(\GR_1 \times \GR_2)$ can be  
identified with a natural, dense subalgebra of the algebra of  
$q$-suspended pseudodifferential operators ``on'' $M$, introduced 
by Melrose.  
\end{example}

\begin{example}\   
\label{ex16}    
If $\GR \to B$ is a bundle of Lie groups, then $\tPS{m}$ consists of  
smooth families of invariant, properly supported, pseudodifferential  
operators on the fibers of $\GR \to B$. For $\GR = B \times G$, a  
trivial bundle of Lie groups, $\GR$ is the product (as groupoids) of a  
smooth manifold $B$, as in Example \ref{ex.manifold}, and a Lie group  
$G$, as in Example \ref{ex.Lie}.  A very important particular case of  
this construction is when $\GR \to B$ is a vector bundle, with the  
induced fiberwise operations. We shall use this example below several  
times.  
\end{example}

\begin{example}\ Again, the only general thing that can be said about fibered   
products is that $\Psi^{-\infty}(\GR_1) \otimes_{\CI(B)}  
\Psi^{-\infty}(\GR_2)$ identifies with a dense subset of  
$\Psi^{-\infty}(\GR_1 \times_B \GR_2)$.  When $\GR_1 = M \times_B M$  
is a fibered pair groupoid and $\GR_2$ is a bundle of Lie groups on  
$B$, then $\Psi^m(\GR_1 \times \GR_2)$ is an algebra considered in  
\cite{NistorIndFam}, and consists of smooth families of  
pseudodifferential operators on $M \times_B \GR_2$ invariant with  
respect to the bundle of Lie groups $\GR_2$.  
\end{example}

\begin{example}[Connes]\ If $\GR$ is the holonomy groupoid associated to    
the foliated manifold $(M,F)$, then $\Psi^*(\GR)$ is the algebra of   
pseudodifferential operators along the leaves of $(M,F)$, considered   
first by Connes \cite{ConnesF}. In fact, our algebra is a little smaller   
than Connes' who considered   
families that are only continuous in the transverse direction.   
These algebras, however, have the same formal properties as our families.    
\end{example}

\begin{example}\ Let $\GR$ be the fundamental groupoid of a    
compact smooth manifold $M$ with fundamental group   
$\pi_1(M)=\Gamma$. If $P=(P_x)_{x \in M} \in \tPS{m}$, then each   
$P_x$, $x \in M$, is a pseudodifferential operator on   
$\widetilde{M}$. The invariance condition applied to the elements $g$   
such that $x=d(g)=r(g)$ implies that each operator $P_x$ is invariant   
with respect to the action of $\Gamma$. This means that we can   
identify $P_x$ with an operator on $\widetilde M$ and that the   
resulting operator does not depend on the identification of $\GR_x$   
with $\widetilde M$. Then the invariance condition applied to an   
arbitrary arrow $g \in \GR$ gives that all operators $P_x$ acting on   
$\widetilde M$ coincide. We obtain $\tPS{m} \simeq   
\Psi^m_{\prop}(\widetilde{M})^\Gamma$, the algebra of properly   
supported $\Gamma$-invariant pseudodifferential operators on the   
universal covering $\widetilde{M}$ of $M$. An alternative definition   
of this algebra using crossed products is given in \cite{Nistor4}.   
\end{example}

\begin{example}\   
\label{exad}  
If $\GR_{\adb}$ is the adiabatic groupoid associated to   
a groupoid $\GR$, then an operator $P \in \Psi^m(\GR_{\adb})$ consists   
of a family $P=(P_{t,x})$, $t \ge 0$, $x \in M$, ($M$ is the space of   
units of $\GR$), such that if we denote by $P_t$ the family   
$(P_{t,x})$, for a fixed $t$, then $P_t \in \tPS{m}$ for $t > 0$ and   
$P_{t}$ depends smoothly on $t$ in this range.  For $t = 0$, $P_0\in   
\Psi^{m}(A(\GR))$, is a family of operators on the fibers of $A(\GR)   
\to M$, translation invariant with respect to the variable in each   
fiber.  Thus, $\Psi^m(A(\GR))$ is one of the algebras    
appearing in Example \ref{ex16}.   
 
In a certain sense $P_t \to P_0$, as $t \to 0$, but this is difficult to   
make precise without considering the adiabatic groupoid. (Actually,   
making precise the fact that the family $P_t$ is smooth at $0$ also is   
precisely the {\em raison d'\^etre} for the algebra of   
pseudodifferential operators on the adiabatic groupoid.)   
   
The best way to formalize this continuity property is the following.   
Consider the evaluation morphisms $e_t : \Psi^m(\GR_{\adb}) \to   
\Psi^m(\GR)$, if $t >0$, and $e_0 : \Psi^m(\GR_{\adb}) \to   
\Psi^m(A(\GR))$. (These morphisms are particular instances of the restriction   
morphisms defined in Equation \eqref{eq.indicial}).    
If $P \in \Psi^0(\GR_{\adb})$, then $\| e_t(P) \|$   
and $\|e_t(P)\|_r$ are continuous in $t$. This was proved by Landsman    
and Ramazan, see \cite{Landsman,LandsmanRamazan,Ramazan}.   
 
Some typical operators in $\Psi^m(\GR_{\adb})$ are obtained by 
rescaling the symbol of a differential operator $D$ on $\GR$. To see 
how this works, note first that there exists a polynomial symbol $a$ on 
$A^*(\GR)$ such that $q(a) = D$, where $q : {\mathcal S}^m(A^*(\GR)) 
\to \tPS{m}$ is the quantization map considered in \cite{NWX}. Let 
$a_t$ be the symbol $a_t(\xi) = a(t\xi)$, for $t > 0$, and also let 
$q_{\adb} : {\mathcal S}^m(A^*(\GR_{adb})) \to \Psi^m(\GR_{\adb})$ be 
the quantization map for the adiabatic groupoid. We can extend $a$ to 
a symbol on $A^*(\GR_{\adb})$ constant in $t$, then $e_t(q_{\adb}(a)) 
= q(a_t)$, if $t > 0$. For $t = 0$, we obtain that $e_0(q_{\adb}(a))$ 
is isomorphic to the operator of multiplication by $a$, after taking 
the Fourier transform along the fibers of $A^*(\GR)$. 
\end{example}

An important class of examples is obtained by integrating suitable Lie 
algebras of vector fields on a manifold $M$ with corners.  This is 
related to Melrose's approach to a pseudodifferential analysis on 
manifolds with corners \cite{MelroseScattering}, though our techniques 
are different in the end.  We thus start with a Lie subalgebra $\cV$ 
of the Lie algebra of all vector fields that are tangent to each 
boundary hyperface of a given manifold $M$ with corners.  The Lie 
algebra $\cV$ can be thought of as determining the degeneracies of our 
operators near the boundary. If $\cV$ is in addition a projective 
$\CI(M)$-module, then, by the Serre-Swan theorem, there is a smooth 
vector bundle $A={}^{\cV}TM\rightarrow M$ together with a smooth map 
of vector bundles $q:A\longrightarrow TM$ such that 
$\cV=q(\Gamma(A))$.  (This will be discussed in more detail in a 
forthcoming book of Melrose on manifolds with corners.) 
  
The next step is to integrate this Lie algebroid $A$, that is, to find  
a Lie groupoid $\GR$ with Lie algebroid $A$. Here, we can follow the  
general method used in \cite{NistorInt}. The integration procedure  
consists in fact of two steps.  Let us denote by $A_{S}$ the  
restriction of $A$ to each open boundary face $S$ of $M$ of positive  
codimension, suppose that we can find differentiable groupoids $\GR_S$  
integrating $A_S$, and let $\GR=\bigcup\GR_S$.  By \cite{NistorInt},  
there exists at most one smooth structure on $\GR$ compatible with the  
groupoid operations.  Whenever such a smooth structure exists, the  
resulting groupoid satisfies $A(\GR) = A$.  Moreover, if the $\GR_S$  
are maximal among all $d$-connected groupoids integrating $A_S$, then  
there is a natural differentiable structure on $\GR$ making it into a  
differentiable groupoid with Lie algebroid $A$. Note that this choice  
for $\GR$ will almost always lead us to non-Hausdorff groupoids and to  
problems related to the analysis on these spaces. Moreover, the vector  
representation will not be injective, in general. The reason is that  
the maximal $d$-connected groupoid integrating a given Lie algebroid  
is much to big. For instance, for the Lie algebroid $TM \to M$, the  
maximal $d$-connected groupoid integrating it is the path groupoid  
\cite{NWX}, not the pair groupoid as expected and usually desired  
\cite{Mackenzie1}.  In particular cases, however, the given Lie  
algebroid $A$ can be integrated directly to a Hausdorff  
differentiable groupoid. These remarks apply to the following two  
examples. These two examples are essentially due to Melrose  
\cite{MelroseScattering} and, respectively, to Mazzeo \cite{maz91}.  A  
groupoid for a special case of Example \ref{excn} (b-calculus on 
manifolds with corners) was constructed in \cite{Monthubert}, see also 
\cite{NWX}.

\begin{example}\    
\label{excn} {\em The ``very small'' $c_n$-calculus.}   
Let $M$ be a compact manifold with corners, and associate to each  
hypersurface $H \subset M$ an integer $c_H \ge 1$.  We also fix a  
defining function for each hypersurface.  Choose also on $M$ a metric  
$h$ such that each point $p \in F$, belonging to the interior of a  
face $F\subseteq M$ of codimension $k$, has a neighborhood $V_p \cong  
V'_p \times [0,\varepsilon)^k$, with the following two properties:\  
the defining function $x_j$ is obtained as the projection onto the  
$j$th component of $[0,\varepsilon)^k$ and the metric $h$ can be  
written as $h = h_F + (dx_1)^2 + \ldots (dx_k)^2$, with $x_1, \ldots,  
x_k$ being the defining functions of $F$ and $h_F$ being a two-tensor  
that does not depend on $x_1,\ldots, x_k$ and restricts to a metric on  
$F$.

Then, we consider on $M$ the vector fields $X$ that in a neighborhood   
of each point $p$, as above, are of the form   
$$   
        X = X_F + \sum_{j=1}^k x_j^{c_j} \pa_{x_j},   
$$   
with $c_j$ being the integer associated to the hyperface $\{x_j =   
0\}$ and $X_F$ being the lift of a vector field on $F$.   
The set of all vector fields with these properties forms a Lie   
subalgebra of the algebra of all vector fields on $M$. We denote this   
subalgebra by ${\mathcal A}(M,c)$. By the Serre-Swan    
theorem, there exists a vector bundle   
$A(M,c)$ such that $\mathcal A(M,c)$ identifies with the space of   
smooth sections of $A(M,c)$.   
   
We want to integrate $A(M,c)$, and to this end, we shall use the approach    
from \cite{NistorInt}.  Let $S = int(F)$ be the   
interior of a face $F \subset M$ of codimension $k$. The restriction    
of $A(M,c)$ to each open face $S$ is then $TS \times \RR^k$, and hence it   
is integrable; a groupoid integrating    
this restriction being, for example $\GR_S = S \times S \times \RR^k$, if    
$F = \overline{S}$ has codimension $k$.    
     
Define then   
$$   
        \GR := \bigcup_F S \times S \times \RR^k,   
$$   
which is a groupoid with the obvious induced structural maps.  As a  
set, $\GR$ does not depend on $c$. Because the groupoids $\GR_S$ are  
not $d$-connected, in general, we cannot use the result of  
\cite{NistorInt} to prove that it has a natural smooth structure, so  
we have to construct this smooth structure directly.  
   
The results of \cite{NistorInt} say that if there exists a smooth  
structure on $\GR$ compatible with its groupoid structure, then it  
must be obtained using certain coordinate charts defined using the  
exponential map.  In our case, the exponential map amounts to the  
following.  
   
Let $\psi_l : (0,\infty) \to \RR$ be $\psi_l(x) = \ln x$, if $l = 1$  
and $\psi_l(x) = x + x^{1-l}/(1-l)$, if $l > 1$. Also, let $\phi_l :  
\RR \times [0,\infty) \to [0,\infty)$ be defined by $\phi_l(t,0) = 0$  
and $\phi_l(t,x) = \psi_l^{-1}( \psi_l(x) + t)$.  In particular,  
$\phi_1(t,x) = e^tx$. Then $\phi_l$ defines a differentiable action of  
$\RR$ on $[0,\infty)$, which hence makes $\RR \times [0,\infty)$ a  
differentiable groupoid denoted ${\mathcal F}_l$.  The Lie algebroid  
of ${\mathcal F}_l$ is generated as a $\CI([0,\infty))$-projective  
module by the infinitesimal generator $\pa_t$ of the action of $\RR$;  
note that the action of $\partial_{t}$ on $\CI([0,\infty))$ under the  
anchor map is given by $f(x)x^{l}\partial_{x}$ for some nowhere  
vanishing (bounded) smooth function $f$.  Consequently, $A({\mathcal  
F}_l)$ is the projective $\CI([0,\infty))$-module generated by $x^l  
\pa_x$.  
   
Assume now that $M = [0,\infty)$ and fix $l\in \NN$. Then ${\mathcal  
F}_l$ is a smooth groupoid integrating $A(M,l)$, by the above remarks.  
Consequently, if $M = [0,\infty)^n$ and $c=(c_1,c_2,\ldots,c_n)$, then  
$\GR := {\mathcal F}_{c_1} \times {\mathcal F}_{c_2}\times \ldots  
\times {\mathcal F}_{c_n}$ satisfies $A(\GR) = A(M,c)$. To integrate  
general Lie algebroids of the form $A(M,c)$ we localize this  
construction. This then gives the following smooth structure on $\GR  
:= \bigcup_S \GR_S$.  
   
We now discuss the general case of a manifold with corners. Locally,  
the smooth structure on $\GR$ is given by the discussion above. Since  
this smooth structure is important in applications, let us try to make  
it more explicit.  Thus, fix an arbitrary point $(p,q,\xi) \in S  
\times S \times \RR^k$, which we want to include in a coordinate  
system. By definition, $p,q \in S$. Choose now a small coordinate  
neighborhood $V_p \cong V'_p \times [0,\varepsilon)^k$ of $p$, with  
$V'_p$ a small open neighborhood of $p \in S$, as above. Choose  
$V_q\cong V'_q \times [0,\varepsilon)^k$ similarly.  We write  
$$   
        z = (z', x_1(z), x_2(z), \ldots, x_k(z))   
$$    
for any $z \in V_p \cap V_q$; this is possible since we can assume  
that $V_p = V_q$ if $p = q$ or that $V_p \cap V_q = \emptyset$ if $p  
\not = q$. Fix $R > 2\|\xi\|$ and choose $\delta>0$ so small that  
$|\phi_l(t, x)| < \varepsilon$ if $|t | \le R$, $x \le \delta$, and $l  
= c_j$, for $j = 1, 2, \ldots, k$.  Here $c_j$ is the constant  
associated to the hyperface $\{x_j = 0\}$.  Then we define a map  
\begin{eqnarray}   
        F : V'_p \times [0,\delta)^k   \times V'_q \times    
        \{ \| \xi \| < R \} &\longrightarrow& \GR: \\ \nonumber   
        (z', y, z'', \xi) &\longmapsto&   
        (z', y, z'', \Phi(\xi,y), p_y(\xi)) \in \GR_{S'}   
\end{eqnarray}   
as follows. Let $y = (y_1, \ldots, y_k)$, $B \subset \{ 1, 2, \ldots,  
k \}$ be the subset of those indices $j$ such that $y_j = 0$, and let  
$p_y : \RR^k \to \RR^B$ be the corresponding projection. The vector  
space $\RR^B$ identifies naturally with the fiber at $(z',y)$ of the  
normal bundle to the open face containing $(z',y)$ (this open face was  
denoted above by $S'=S'(z',y)$). For $y = (y_1, y_2, \ldots, y_k)$ and  
$\xi = (\xi_1, \xi_2, \ldots, \xi_k)$, the map $\Phi$ is then given by  
$$   
        \Phi(\xi, y) = ( \phi_{c_1}(\xi_1, y_1), \phi_{c_2}(\xi_2, y_2),    
        \ldots, \phi_{c_k}(\xi_k, y_k))\,.   
$$   
   
We shall denote by $\GR(M,c)$ the smooth groupoid constructed above.   
    
Fix a face $F \subset M$ of codimension $k$. By construction, $F$ is  
an invariant subset of $M$ and hence we can consider the restriction  
maps $\inn_F$ defined in Equation \eqref{eq.indicial}. The range of  
these restriction (or indicial) maps is in related to the algebras  
$\Psi^\infty(\GR(M,c))$.  The precise relation is the following.  
   
Each hyperface $H'$ of $F$ is a connected component of $H \cap F$, for  
a unique hyperface $H$ of $M$. Then we associate to $H'$ the integer  
$c_H \ge 1$.  We denote by $c'$ the collection of integers obtained in  
this way.  Then the restriction of $A(M,c)$ to $F$ is isomorphic to  
$A(F,c') \times \RR^k$.  {}From this we obtain that $\GR(M,c)\vert_{F}  
\cong \GR(F,c') \times \RR^k$. The restriction maps thus become  
$$   
        \inn_F : \Psi^m(\GR(M,c)) \to \Psi^m(\GR(F,c') \times \RR^k).  
$$    
The right hand side algebras are closely related to the ``$k$-fold  
suspended algebras'' of Melrose.  
   
The analytic properties of the algebras $\Psi^\infty(\GR(M,c))$ will  
be studied again in Section \ref{Sec.Examples.III}.  
\end{example}

\begin{example}\label{ex.boundary.f1}   
Let $M$ be a compact manifold whose boundary $\partial M$ is the total  
space of a locally trivial fibration $p:\partial M \longrightarrow B$  
of compact smooth manifolds.  A smooth vector field on $M$ is called  
an {\em edge vector field} if it is tangent to the fibers of $p$ at  
the boundary. The Lie algebra $\VeM$ of all edge vector fields is a  
projective $\CI(M)$-module, and hence, by the Serre-Swan theorem  
\cite{Karoubi}, it can be identified with the space of all $\CI$  
sections of a smooth vector bundle $\eTM\rightarrow M$ that comes  
equipped with a natural map $\eTM\longrightarrow TM$ \cite{maz91}  
making $A:=\eTM$ into a Lie algebroid. A pseudodifferential calculus  
adapted to this setting was constructed by Mazzeo \cite{maz91}, and,  
in a slightly different way by Schulze \cite{braun}. To integrate  
$\eTM$, we shall use the methods of \cite{NistorInt}. 
   
Let $M_0 := M \smallsetminus \pa M$ and notice that $A\vert_{M_0}  
\cong TM_0$. We can integrate this restriction to the pair groupoid:\  
$\GR_{M_0} := M_0 \times M_0$. The restriction of $A$ to the boundary  
is the crossed product of another Lie algebroid with $\RR$:\  
$$   
        A\vert_{\pa M} \cong (T_{vert}\pa M \times  TB) \rtimes \RR.   
$$   
It is worthwhile do describe this restriction more precisely. As a    
vector bundle, $A$ is the direct sum of three vector bundles:\   
$T_{vert}\pa M$ (the vertical tangent bundle   
to the fibers of $\pa M \to B$),\ $p^*(TB)$  (the pull-back of the   
tangent bundle of $B$), and a trivial, one-dimensional real vector bundle.   
Thus, every section of $A$ can be represented as a   
triple $(X,Y,f)$, where $X$ is a vector field on $\pa M$, tangent to   
the fibers of $\pa M \to B$, $Y$ is a section of $p^*(TB)$, which is   
convenient to be represented as a section of the quotient    
$T \pa M /T_{vert} \pa M$, and $f \in \CI(\pa M)$. Let $\nabla$ be the   
Bott connection on $p^*(TB)$. The Lie algebra structure on $\Gamma(A)$    
is then   
$$   
        [(X,Y,f), (X_1,Y_1,f_1)] = ([X,X_1], \nabla_X(Y_1) + fY_1    
        - \nabla_{X_1}(Y) - f_1Y, 0).    
$$   
   
Let $G \to B$ be the bundle of Lie groups obtained as the  
cross-product of the bundle of commutative Lie groups $TB$ with $\RR$,  
the action of $t \in \RR$ being as multiplication with $e^t$. The Lie  
algebroid (or the bundle of Lie groups associated to this bundle of  
Lie groups) is $A(G) = TB \oplus \RR$, with the bracked defined as  
above:\ $[(Y,f), (Y_1,f_1)] = (fY_1 - f_1Y, 0)$.  Then we can write  
$A\vert_{\pa M} = T_{vert}\pa M \times_B A(G)$. This writing  
immediately leads to a groupoid integrating $A\vert_{\pa M}$, namely  
the fibered product of a groupoid integrating $T_{vert}\pa M$ and a  
groupoid integrating $A(G)$. We can choose these groupoids to be the  
fibered pair groupoid $\GR_1 : = \pa M \times_B \pa M$ and,  
respectively, $G$. The resulting groupoid integrating $A\vert_{\pa M}$  
is then $\GR_{\pa M} := \GR_1 \times_B G \to B$, invariant with  
respect to the action of $G$ by right translations.  The resulting  
algebra of pseudodifferential operators will be an algebra of smooth  
families acting on the fibers of $\pa M \times_B G$, invariant with  
respect to $G$.  
   
To obtain a groupoid integrating $A$, it is enough to show that the  
disjoint union $\GR := \GR_1 \cup (M_0 \times M_0)$ has a smooth  
structure compatible with the groupoid structure. This smooth  
structure is obtained using the following coordinate charts. Let $x$  
be a boundary defining function on $M$, and fix $q \in \pa M$ and a  
neighborhood $V_q \cong V'_q \times [0,\varepsilon)$ such that the  
defining functions $x$ becomes the second projection on $V_q$ and  
$V'_q$ is a neighborhood of $q$ in $\pa M$.  We replace $V_q$ with a  
smaller neighborhood, if necessary, so that there exists a fiber  
preserving diffeomorphism $\phi : B_1 \times B_2 \to V'_q$, $\phi(0,0)  
= q$, from a product of two small open balls in some Euclidean spaces  
(so $p$ becomes the first projection with respect to the  
diffeomorphism $\phi$). Let $q'\in \pa M$ be a second point, chosen  
such that $p(q') = p (q)$, and choose a diffeomorphism $\phi' : B_1  
\times B'_2 \to V_{q'}$ as above.  We can assume that $p \circ \phi =  
p \circ \phi'$ and $\phi = \phi'$, if $q = q'$, or that $V_q$ and  
$V_{q'}$ are disjoint. Then $\phi$ and $\phi'$ define a diffeomorphism  
$$   
        \phi \times_B \phi' : B_1 \times B_2 \times B'_2 \to    
        \pa M \times_B \pa M,   
$$    
explicitly,    
$$   
        \phi \times_B \phi'(b_1, b_2, b'_2) = (\phi(b_1,b_2),  
      \phi'(b_1,b'_2)) \in \pa M \times_B \pa M \subset \pa M \times  
      \pa M.  
$$     
   
We identify $B_1$ with the fiber of $TB$ at $p(q) = p(q')$ such that   
$p(q)$ corresponds to $0$, and  we let   
$$   
        \Phi : B_1 \times B_2 \times [0,\delta) \times B_1 \times B_2 \times (-R,R)   
        \to \GR   
$$   
be given by   
$$   
        \Phi(b_1,b_2, 0, b'_1, b'_2, t) = \big( \phi \times_B  
        \phi'(b_1, b_2, b'_2), b'_1,t) \in \pa M \times_B \pa M \times  
        T_{p(q)}B \times \RR \in \GR_1  
$$   
or by   
\begin{eqnarray*}  
        \Phi(b_1, b_2, s, b'_1, b'_2, t) &=& (\phi(b_1,b_2), s,  
        \phi'(b_1 + sb'_1,b'_2), se^t) \\ &\in& V'_{q} \times  
        (0,\varepsilon) \times V'_{q'} \times (0,\varepsilon) \subset  
        M_0 \times M_{0}.  
\end{eqnarray*}  
   
The restriction at the boundary map $\inn_{\pa M}$ defined in Equation  
\eqref{eq.indicial} becomes a map  
$$   
        \inn_{\pa M} : \tPS{\infty} \to \Psi^\infty(\GR_1).   
$$   
The range of this map consists of families of pseudodifferential  
operators that act on the fibers of $\pa M \times_B G \to B$ and are  
$G$-invariant with respect to the action of $G$ by right translations.  
\end{example}   
 
We expect that the above example will be useful for the question from 
\cite{FreedWitten} on the Bojarsky additivity formula for the real 
index of families of elliptic operators. (See also Nicolaescu's paper 
\cite{Nicolaescu}.) Also, it will probably be useful for a certain 
approach in the study of the $S^1$-equivariant Dirac operators  
\cite{NistorDIR} and \cite{Lott}.

\section{Geometric operators\label{Sec.GO}}

For two vector bundles $E_{0}, E_{1}$ on $M$, we shall denote by  
$\Diff{\GR;E_0,E_1}$ the space of differential operators $D :  
\Gamma(\GR;r^{*}E_0) \to \Gamma(\GR;r^{*}E_1)$ with smooth  
coefficients that differentiate only along the fibers of $d: \GR \to  
M$, and which are right invariant. Thus, $\Diff{\GR;E_0,E_1}$ is  
exactly the space of differential operators in  
$\Psi^m(\GR;E_0,E_1)$. The elements of  
$\Diff{\GR;E_0,E_1}$ will be called {\em differential operators on $\GR$}.   
   
In this section, we define and study the geometric differential  
operators on a given differentiable groupoid $\GR$.  For the  
definition of most of these operators, we shall need a metric on  
$A:=A(\GR)$.  
   
To define the de Rham operator, however, we need no metric.  Denote  
then by $\lambda_q=\Lambda^q T^*_{vert} \GR$ the $q$th exterior power  
of the dual of $T_{vert}\GR$, the vertical tangent bundle to the  
fibers of $d : \GR \to M$.  Recall that $\GR_x$ denotes $d^{-1}(x)$  
throughout this paper. Then the de Rham differential $d :  
\Gamma(\GR_x,\lambda_q) \to \Gamma(\GR_x,\lambda_{q+1})$ is invariant  
with respect to right translations, and hence it defines an operator  
$d \in \bigoplus_q \Diff{\GR;\Lambda^qA^*, \Lambda^{q+1}A^*} \subset  
\Diff{\GR; \Lambda^*A^*}$.  
   
Let as before $\pi$ be the representation $\pi : \Psi^m(\GR;E_0,E_1)  
\to \Hom(\Gamma(E_0),\Gamma(E_1))$ given by the formula $\pi(P)(f)  
\circ r = P(f \circ r)$.  (Recall that we called this representation  
the {\em vector representation}.) The complex determined by the  
operators $\pi(d)$ then computes the {\em Lie algebroid cohomology} of  
$A$, by definition, see \cite{Mackenzie1}.  
    
We shall use later on the explicit form of $d$ for the adiabatic  
groupoid associated to $\GR$. Recall that an operator $P$  
(differential or pseudodifferential) on the adiabatic groupoid  
$\GR_{\adb}$ associated to $\GR$ consists of a family $P=(P_t)$, such  
that, in particular, $P_t \in \tPS{m}$, for $t >0$. The de Rham  
operator $d^{\GR_{\adb}}=(d_t)$ is such that $d_t = td$, for $t > 0$.  
   
Of course, more operators are obtained if we consider a metric on  
$A:=A(\GR)$. Here, by ``metric on $A$'' we mean a positive definite  
bilinear form on $A$, as usual. The metric on $A$ then makes each  
$\GR_x$ a Riemannian manifold, naturally, due to the isomorphisms  
$T\GR_x \cong r^*(A)$ as vector bundles on $\GR_x$.  Moreover, right  
translation by an element of $\GR$ is an isometric  
isomorphism. Because of this, every geometric differential operator  
associated naturally to a Riemannian metric will be invariant with  
respect to the right translation by an element of $\GR$, and hence  
will define an element in $\Diff{\GR;E_0,E_1}$, for suitable vector  
bundles $E_i$. We shall not try to formulate this in the greatest  
generality, but we shall apply this observation to particular  
operators that appear more often in practice.  
   
For example, the metric allows us to define the Hodge $*$-operator,  
which then leads to the signature operator $d \pm *d* \in  
\Diff{\GR,\Lambda^*A}$.  Also, the metric gives rise to an inner  
product on $\Lambda^*A$ and hence to an adjoint to $d$, denoted $d^*$,  
which then in turn allows us to define the Euler operator $d + d^* \in  
\Diff{\GR;\Lambda^*A^*}$. Similarly, one obtains the Hodge Laplacians  
$\Delta_p \in \Diff{\GR;\Lambda^p A^*}$, as components of the square  
of the Euler operator $d + d^*$. We write $\Delta_p^\GR$, $d^\GR$,  
... for these operators when we want to stress their dependence on  
$\GR$.  
   
We now turn to Dirac and generalized Dirac operators. This requires us  
to introduce the (generalization to groupoids of the) Levi-Civita  
connection.  
   
For $X \in \Gamma(A)$, we shall denote by $\tilde X$ its lift to a  
right invariant, $d$-vertical vector field on $\GR$. Let $\nabla^x :  
\Gamma(T_{vert}\GR_x) \to \Gamma(T_{vert}\GR_x \otimes  
T_{vert}^*\GR_x)$ be the Levi-Civita connection associated to the  
induced metric on $\GR_x$. Then for any $X \in \Gamma(A)$, we obtain a  
smooth, right invariant family of differential operators  
$$   
        \nabla^x_{\tilde{X}} : \Gamma(T_{vert}\GR_x) \to  
        \Gamma(T_{vert}\GR_x).  
$$    
We denote the induced differential operator in $\Diff{\GR, A}$ simply  
by $\nabla_X$. For all smooth sections $X$ and $Y$ of $A$, there  
exists another smooth section $Z$ of $A$ such that $\nabla_X (\tilde  
Y) = \tilde Z$.  
   
Suppose now that $A$ is $spin$, that is, that $A$ is orientable and  
the bundle of orientable frames of $A$ lifts to a principal $Spin(k)$  
bundle ($k$ being the rank of $A$). Suppose $k=2l$ is even, for  
simplicity, and let $S = S_+ \oplus S_-$ be the spin bundle associated  
to the given spin structure and the spin representation of  
$Spin(k)$. As in the classical case, the Levi-Civita connection on the  
frame bundle of $r^*(A)$ lifts to a connection $\nabla^S$ on  
$S$. Moreover, this connection involves no choices (it is uniquely  
determined by the spin structure), and hence $\nabla^S$ is right  
invariant, in the obvious sense. Thus, if $X$ is a section of $A$ and  
$\tilde X$ is its lift to a right invariant, $d$-vertical vector field  
on $\GR$, then $\nabla^S_{\tilde X}$ is a right invariant differential  
operator, and hence it is in $\Diff{\GR, S}$. We denote by $\Dir^S$  
the induced Dirac operator on the spaces $\GR_x$, which will then form  
a right invariant family, and hence $\Dir^S \in \Diff{\GR;S}$. (We  
shall write $\Dir^S_{\GR}$ on the few occasions when we shall need to  
point out the dependence of this operator on the groupoid $\GR$.)  
   
Let $\Cliff(A)$ be the bundle of Clifford algebras associated to $A$  
and its metric.  We shall use the metric to identify $A^*$ with $A$,  
so that $\Cliff(A^*)$ becomes identified with $\Cliff(A)$. The same  
construction as above then applies to a $\Cliff(A)$-module $W$ endowed  
with a right invariant, admissible connection $\nabla^W$ (see below)  
on each of its restrictions to $\GR_x$. Denote by $c : \Cliff(A) \to  
End(W)$ the Clifford module structure on $W$. Because $A \subset  
\Cliff(A)$, we also obtain a bundle morphism $A \to End(W)$ still  
denoted $c$. Recall then that $\nabla^W$ is an {\em admissible  
connection} if, and only if,  
$$   
        \nabla_X^W(c(Y)\xi) = c(\nabla_XY)\xi + c(Y)\nabla^W_X(\xi),   
$$   
for all $\xi \in \Gamma(r^*(W))$ and all $X,Y \in \Gamma(r^*(A))$, the   
second connection being the Levi-Civita connection discussed   
above. Then we obtain as in the classical case a Dirac operator   
$\Dir^W_x$ on $\GR_x$, acting on sections of $r^*(W)$. The right   
invariance of the connection $\nabla^W$ guarantees that   
the family $\Dir^W_x$ is   
right invariant, and hence that it defines an element in $\Diff{\GR;W}$.   
   
It is a little bit trickier to define the generalized Dirac operator  
associated to a $\Cliff(A)$-module $W$, if no admissible connection is  
specified on $W$. This is because it is not clear a priori that right  
invariant admissible connections exist at all. Our next goal then is  
to prove that this is always the case, as it is for Clifford modules  
on Riemannian manifolds.  
   
We shall work with complex $\Cliff(A)$-modules, for simplicity. Also,  
we assume that $A$ is even dimensional.  Cover $M$ with contractible  
open sets $U_{\alpha}$. Then $A\vert_{U_\alpha}$ has a trivialization  
$A\vert_{U_\alpha} \simeq U_{\alpha} \times \RR^{2l}$, which we can  
assume to preserve the metric. Then $\Cliff(A)\vert_{U_{\alpha}}  
\simeq U_{\alpha} \times M_{2^l}(\CC)$ and $W\vert_{U_\alpha} \simeq  
U_{\alpha} \times V \simeq \CC^{2^l} \otimes V_0$, with $V_0$ an  
additional vector bundle, which is acted upon trivially by the  
Clifford algebra, and hence only serves to encode the local  
``multiplicity'' of the $\Cliff(A)$-module $W$. As in the classical  
case, we first define the admissible connection locally, using the  
above trivialization, and then we glue them using a partition of  
unity.  However, in our groupoid setting we need to work a little bit  
more to make sense of what the ``local definition'' means. More  
precisely, all definitions will be given not on $U_\alpha$ itself, but  
on $r^{-1}(U_\alpha)$. Once we realize this, everything carries over  
from the case of a Riemannian manifold to that of a differentiable  
groupoid. For completeness, we now review this construction in our  
case.  
   
The trivialization of $U_\alpha$ gives an orthonormal family of  
sections $X_1, X_2,\ldots,X_k$ of $A$ over $U_\alpha$. Then, we obtain  
smooth functions $\Gamma^a_{bc}$ on $U_\alpha$ such that, working {\em  
always} over $r^{-1}(U_\alpha)$,  
$$   
        \nabla_{\tilde X_i} \tilde X_j = \sum_h   
        (\Gamma^h_{ij} \circ r) \tilde X_h.   
$$   
(Compare with \cite{LM}.)  Fix a basis $(e_t)$, $t = 1, \ldots, 2^lm$,  
of $V$, where $V$ is the vector space appearing in the isomorphism  
$W\vert_{U_\alpha} \simeq U_\alpha \times V$. We shall denote by  
$\tilde e_t := e_t \circ r$ the induced basis of $r^*(W)$ on  
$r^{-1}(U_\alpha)$.  The point of these choices is, of course, that  
the matrix of the multiplication operator $c(\tilde X_j)$ in the basis  
$\tilde e_t$ consists of {\em constant} functions.  Using the  
functions $\Gamma^h_{ij}$ and the Clifford multiplication map $c: A  
\to End(W)$, we define a connection $\nabla^{x,W,\alpha}$ on the  
restriction of $r^*(W)$ to $\GR_x \cap r^{-1}(U_\alpha)$ by the  
formula  
\begin{equation}   
        \nabla_{\tilde X_h}^{x,W,\alpha} \tilde e_t := \frac{1}{4}\sum_{a,b}   
        (\Gamma_{ha}^b \circ r) c(\tilde X_a)c(\tilde X_b) \tilde e_t.   
\end{equation}   
    
Let $\phi_\alpha \in \CI(M)$ be a $\CI$-partition of unity subordinate  
to the covering $U_\alpha$.  Then $\tilde \phi_\alpha :=  
\phi_\alpha\circ r$ is a partition of unity subordinate to  
$r^{-1}(U_\alpha)$.  We define a connection $\nabla^{x,W}$ on the  
restriction of $W$ to $\GR_x$ by the formula  
$$   
        \nabla^{x,W}_{\tilde X}(\xi) = \sum_{\alpha}   
        \nabla^{x,W,\alpha}(\tilde X) (\tilde \phi_\alpha \xi).   
$$   
By the definition, $\nabla^{x,W}$ is an admissible connection on the  
restriction of $r^*(W)$ to $\GR_x$.

\begin{proposition}\    
Let $W \to M$ be a complex vector bundle that is a $\Cliff(A)$-module.  
Then we can find an admissible connection $\nabla^{x,W}$ on the  
restriction of $r^*(W)$ to $\GR_x$, for any $x \in M$, such that for  
each $X \in \Gamma(A)$, the operators $\nabla^{x,W}_{\tilde X}$ form a  
smooth, $\GR$-invariant family of differential operators on $r^*(W)$,  
and hence they define an element $\nabla_X^W$ in $\Diff{\GR;W}$. If  
$S$ is a spin bundle, then we can take this connection to be the  
Levi-Civita connection.  
\end{proposition}

\begin{proof} This is just the summary of the above discussion. \end{proof}

It follows from the above proposition that if we consider on each  
$\GR_x$ the Dirac operator determined by the connection  
$\nabla^{x,W}$, then we obtain an invariant family of differential  
operators, which hence defines an operator $\Dir^W_\GR \in  
\Diff{\GR;W}$, the Dirac operators on $\GR$ associated to $W$ and the  
given admissible connection. (When the groupoid $\GR$ is clear from  
the context, we shall drop the subscript $\GR$.)  
   
We can also regard the admissible connection on a $\Cliff(A)$-module  
$W$ as an operator $\nabla^W \in \Diff{\GR;W, W\otimes A^*}$. If we  
denote by $c \in Hom(W\otimes A^*, W)$ the Clifford multiplication,  
then, as in the classical case $\Dir^W = c \circ \nabla^W$. We can  
also generalize the local description of Dirac operators. Let $M =  
\bigcup U_\alpha$ be a covering of $M$ by open subsets which  
trivializes the bundle $A=A(\GR)$, and choose a partition of unity  
$\phi_\alpha^2$ subordinate to $U_\alpha$.  On each $U_\alpha$, we  
choose a local {\em orthonormal} basis $X_1, \ldots, X_k$ of $A$ and  
define $X^\alpha_j = \phi_\alpha X_j$. Then  
\begin{equation}   
        \Dir^W = \sum_{\alpha,j} c(X^\alpha_j) \nabla^W_{X^{\alpha}_j}.   
\end{equation}   
   
As in the classical case of a Riemannian manifold, the space of  
$\GR$-invariant, admissible connections $\nabla^{x,W}$ on $r^*(W)$ is  
an affine space with model vector space the space of skew-adjoint  
elements in the space of $\Cliff(A)$-linear endomorphisms of $W$.  
   
A feature specific to the groupoid case, however, is that all the  
above constructions and operators are compatible with restrictions to  
compact, $\GR$-invariant subsets of $M$. (Recall that a subset $Y  
\subset M$ of the space of units of $\GR$ is $\GR$-invariant if, and  
only if, $d^{-1}(Y) = r^{-1}(Y)$.) For instance, consider a  
$\Cliff(A)$ bundle $W$ on $M$ with admissible connection  
$\nabla$. Then $W$ restricts to a $\Cliff(A\vert_Y)$ module on  
$Y$. {}From this observation we get that the Dirac operator on $\GR$  
associated to the $\Cliff(A)$-module $W$ will restrict to the Dirac  
operator on $\GR_Y := d^{-1}(Y)$ associated to the  
$\Cliff(A\vert_Y)$-module $W\vert_Y$. Formally,  
\begin{equation}   
        \inn_Y(\Dir^W_{\GR}) = \Dir^{W_Y}_{\GR_Y}.   
\end{equation}   
Similarly,   
\begin{equation}   
\label{lapl}   
        \inn_Y(\Delta_p^\GR) = \Delta_p^{\GR_Y}\,, \quad \inn_Y(d^\GR)  
        = d^{\GR_Y},  
\end{equation}   
and so on. This leads, as we shall see in the following sections, to  
Fredholmness criteria for these various operators in terms of the  
invertibility of the corresponding operators associated to proper,  
invariant, closed submanifolds.

\section{Sobolev spaces\label{Sec.Sobolev}}

Throughout this section, we assume for simplicity, that the space $M$  
of units of a given groupoid $\GR$ is compact.  All the definitions  
and results extend to the case of sections of a Hermitian vector  
bundle $E$ and operators acting on sections of $E$. For simplicity,  
however, we shall discuss in detail only the case where $E$ is the  
one-dimensional, trivial bundle.  
  
The notation $E$, sometimes is used in this section to denote the  
identity element of operator algebras, in this  
section.

Consider a bounded, non-degenerate representation  
$\varrho:\tPS{-\infty}\longrightarrow\End(\cH)$.  Theorem  
\ref{Theorem.EXT} then gives a natural extension of $\varrho$ to a  
bounded $*$-representation of $\tPS{\infty}$. Here ``bounded'' refers  
to the fact that the order zero operators act by bounded operators on  
$\mathcal{H}$, see Definition \ref{def.bounded.rep}.  Recall that  
$\varrho$ is non-degenerate if the space  
$\cHu:=\varrho(\tPS{-\infty})\cH$ is dense in $\mathcal{H}$.  The best  
we can hope for formally self-adjoint operators $A=A^{*} \in \tPS{m}$  
is that they are essentially self-adjoint unbounded operators on  
$\mathcal{H}$.  This is in fact the case for elliptic operators; for  
$m>0$, a formally self-adjoint, elliptic operator $A=A^{*}\in\tPS{m}$  
leads to densely defined, essentially self-adjoint operator  
$$   
        \varrho(A):\cHu\longrightarrow\cH.   
$$    
This will allow us to freely use functional calculus for self-adjoint  
operators later on in this section.  Fix a bounded, non-degenerate  
representation $\varrho$ as above.

Note that under the assumptions above, the unit  
$E:=(\id_{\GR_{x}})_{x\in M}$ belongs to $\tPS{0}$ with  
$\varrho(E)=\id_{\cH}$.  Let further $\cHmu:=\cHu^{*}$ be the  
algebraic dual of $\cHu$, and $T:\cH\hookrightarrow\cHmu:h\mapsto  
T_{h}$ be the natural, anti-linear embedding.  As in the classical  
case, $\varrho$ induces a multiplicative morphism  
$\rs:\tPS{\infty}\longrightarrow\End(\cHmu)$ by  
$$   
[\rs(A)u](\xi):=u(\varrho(A^{*})\xi)$$   
for    
$A\in\tPS{\infty}$, $u\in\cHmu$, and $\xi\in\cHu$.   
Chasing  definitions yields   
\begin{eqnarray}\nonumber   
        \rs(A)\circ T= T\circ \varrho(A) &:& \cHu\longrightarrow\cHmu  
        \mbox{ for } A\in\tPS{\infty}\,, \mbox{ and }\\ \rs(A)\circ T=  
        \label{gl2} T\circ \varrho(A) &:&  
        \cH\hspace{2ex}\longrightarrow\cHmu \mbox{ for }  
        A\in\tPS{0}\,.  
\end{eqnarray}   
Recall that an unbounded, closable operator $S:\cH\supseteq  
\cD(S)\longrightarrow\cH$ is called {\em essentially self-adjoint} if  
$\cD(S)$ is dense in $\cH$, and $\overline{S}=S^{*}$ where  
$\overline{S}$ denotes the minimal closed extension of $S$, and  
$S^{*}$ its adjoint in the sense of unbounded operators.  The proof of  
the following proposition is appropriately adapted from \cite[Theorem  
26.2]{shubin}.

\begin{proposition} \label{ess}   
Let $m>0$, and $A=A^{*}\in\tPS{m}$ be elliptic. Then the unbounded  
operator $\varrho(A):\cH\supseteq\cHu\longrightarrow\cH$ is  
essentially self-adjoint. Moreover,  
\begin{equation}\label{domain}   
        \cD\left(\overline{\varrho(A)}\right)=   
        \cD(\varrho(A)^{*})=   
        \left\{h\in\cH:\rs(A)T_{h}\in T\cH\right\}\,.   
\end{equation}   
\end{proposition}

\begin{proof}\ For brevity, let $\cD$ be the space on  
the right-hand side in \eqref{domain}. Also, let  
$h\in\cD(\varrho(A)^{*})$.  Then we get for all $\xi\in\cHu$  
$$   
        \rs(A)T_{h}(\xi)= \left(\varrho(A)\xi,h\right)  
        =\left(\xi,\varrho(A)^{*}h\right)= T_{\varrho(A)^{*}h}(\xi)\,,  
$$   
i.e.\ $\cD(\varrho(A)^{*})\subseteq\cD$, and   
$\rs(A)T_{h}=T_{\varrho(A)^{*}h}$.   
   
On the other hand, for $h\in\cD$, there exists $g\in\cH$ such that   
for all $\xi\in\cHu$    
$$   
        \left(\varrho(A)\xi,h\right)=   
        \rs(A)T_{h}(\xi)=T_{g}(\xi) =(\xi,g)\,,   
$$   
hence, $h\in\cD(\varrho(A)^{*})$ which gives the second equality in  
\eqref{domain}.  By \cite[Theorem 26.1]{shubin}, it remains to show  
$$  
        N(\varrho(A)^{*}\pm i \id_{\cH})\subseteq 
        \cD\left(\overline{\varrho(A)}\right)\,.   
$$  
Because of $m>0$, $A\pm i E\in\tPS{m}$ is elliptic; by the usual 
symbolic argument we get $B_{\pm}\in\tPS{-m}$ satisfying 
$E-B_{\pm}(A\pm iE)=:R_{\pm}\in\tPS{-\infty}$.  Furthermore, for 
$\xi\in N(\varrho(A)^{*}\pm i \id_{\cH})\subseteq\cD(\varrho(A)^{*})$ 
another definition chase yields as before $$ \rs(A\pm iE)T_{\xi}= 
T((\varrho(A)^{*}\mp i \id_{\cH})\xi)=0\,, $$ thus, 
$$  
        T_{\xi}=\rs(R_{\pm})T_{\xi}=T_{\varrho(R_{\pm})\xi}\in T\cHu  
$$ 
because of $R_{\pm}\in\tPS{-\infty}$ and \eqref{gl2}.  Since we have 
$\cHu\subseteq\cD\left(\overline{\varrho(A)}\right)$, this completes 
the proof. 
\end{proof}

Let us now define Sobolev spaces in the setting of groupoids using the  
powers of an arbitrary positive element $D \in \tPS{m}$, $m >0$, as  
customary.  The necessary facts that imply independence of $D$ are  
contained in the following theorem (and the lemmata leading to its  
proof). Also, the following theorem will allow us to reduce certain  
questions about operators of positive order to operators of order  
zero.  
    
We shall write $P \ge 0$ if $P=P^{*} \in \tPS{m}$ is such that  
$(\varrho(P)\xi,\xi) \ge 0$ for all $\xi \in \cHu$ and for every  
non-degenerate representation $\varrho$ of $\tPS{\infty}$ on  
${\mathcal H}$. Also, we shall write $ A \ge B$ if $A - B \ge 0$. For  
$Q \in \alg\GR$, we write $Q = P^{-1}$, if, and only if,  
$\varrho(P)\varrho(Q) = \id_{\cH}$ and $\varrho(Q)\varrho(P) \subseteq  
\id_{\cH}$, for every non-degenerate, bounded representation   
$\varrho$. Then, for $s > 0$, $P^{-s}$ stands for $(P^{-1})^{s}$.

\begin{theorem}\label{theorem.red}\    
Fix a differentiable groupoid $\GR$ whose space of units, $M$, is  
compact. Let $D \in \tPS{m}$, $m > 0$, be such that $D \ge E$ and  
$\sigma_m(D) > 0$.  Then $D^{-s} \in \ideal\GR$, for all $s >  
0$. Moreover, if $P$ has order $\le k$, then $PD^{-k/m} \in \alg\GR$.  
\end{theorem}

The proof will consist of a sequence of lemmata.

\begin{lemma}\label{lemma.3}\ Fix arbitrarily a metric on $A(\GR)$,    
and let $B= E + \Delta$, where $\Delta$ is the positive Laplace  
operator on functions. Then $B$ is invertible in the sense above, and  
we have $B^{-1} \in \ideal\GR$.  
\end{lemma}

\begin{proof}\  
Let $D_{t} = E + t^2 \Delta$, $t>0$. We shall prove first that, for  
small $t$, there exists $Q_t \in \ideal\GR$ such that  
$\varrho(Q_t)\varrho(D_t)\subseteq\id_{\cH}$ and $  
\varrho(D_t)\varrho(Q_t) = \id_{\cH}$, for all non-degenerate  
representations $\varrho$ on $\cH$.  
  
Because the family $(td)$, $t > 0$, extends to a first-order  
differential operator on $\GR_{\adb}$, the adiabatic groupoid of  
$\GR$, we obtain that $t^2\Delta = (td)^* (td)$ induces an element in  
$\Psi^{2}(\GR_{\adb})$, which explains the choice of the power  
$t^2$.  
  
To be precise, let $E\in\tPS{0}$, $E_{\adb}\in\Psi^{0}(\GR_{\adb})$, 
and $E_{0}\in\Psi^{0}(A(\GR))$ be the identity elements.  If $e_t$, $t 
\ge 0$, denotes the evaluation map as in Example \ref{exad} (so that, 
in particular $e_{t}: \Psi^{\infty}(\GR_{\adb}) \to 
\Psi^{\infty}(\GR)$, $t >0$), then we have $e_{t}(E_{\adb})=E$ for 
$t>0$, and $e_{0}(E_{\adb})=E_{0}$.  Thus, the family $(D_t)$ leads to 
an element $D\in\Psi^{2}(\GR_{\adb})$.  Choose a quantization map $q$ 
for $\GR_{\adb}$ as in \cite{NWX}, and denote by $|\xi|$ the metric on 
$A^{*}(\GR)$, so that the principal symbol of $\Delta$ is $|\xi|^2$. 
Then the function $p(t,\xi) := ( 1 + |\xi|^{2})^{-1}$ is an order two 
symbol on $A^*(\GR_{\adb})$, see Example \ref{exad}, 
and $F:=q(p)D\in\Psi^{0}(\GR_{\adb})$ satisfies $e_{0}(F)=E_{0}$. 
>From the results \cite{Landsman,LandsmanRamazan,Ramazan}, we know that 
the function $t\mapsto\|e_{t}(F - E_{\adb})\|$ is continuous at $0$ 
(in fact everywhere, but that is all that is needed), and hence 
$e_{t}(F)$ will be invertible in $\alg{\GR}$ for $t$ small.  We define 
then 
$$ 
        Q_{t}:=e_{t}(F)^{-1}e_{t}(q(p)) \in \alg\GR \tPS{-2}  
        \subseteq \ideal\GR\,,  
$$  
and a straight-forward computation gives $\varrho(Q_t) \varrho(D_t) 
\xi = \xi$, for $\xi\in\cHu$, a dense subspace of $\cH$,  and for $t > 0$ but  
small. Since  $\varrho(D_t)$ is (essentially) self-adjoint, we obtain 
that $\varrho(Q_t)$ is the inverse of (the closure of) 
$\varrho(D_t)$. This means $Q_t = D_t^{-1}$, for $t > 0$ but small, 
according to our conventions. 
  
Let now $h_{t}(y) = (1+t^{2}y)^{-1}$ and $\varepsilon>0$ be  
arbitrary. Then there exists a continuous function  
$g_{\varepsilon,t}:[0,1]\longrightarrow[0,1]$ with $g(0)=0$ such that  
$h_{t}=g_{\varepsilon,t}\circ h_{\varepsilon}$, and we obtain  
$D^{-1}_{t}=g_{\varepsilon,t}(D_{\varepsilon}^{-1}) \in\ideal{\GR}$ by  
the composition property of the functional calculus for continuous  
functions, for $\varepsilon$ small enough.  Because of $B=D_1$ this  
completes the proof.  
\end{proof}

\begin{lemma}\label{lemma.4}\    
Let $D \in \tPS{m}$ be elliptic with $\sigma_m(D) > 0$. Then, for each  
$A \in \tPS{m}$ and for each bounded representation $\varrho$ of  
$\tPS{\infty}$ on $\mathcal H$, we can find $C_A \ge 0$ such that  
$\|\varrho(A)f \| \le C_A(\|f\| + \|\varrho(D) f\|)$, for all $f \in  
\cHu:=\varrho(\tPS{-\infty})\cH$.  
\end{lemma}

\begin{proof}\ The proof is the same as that of the   
boundedness of operators of order zero, using H\"ormander's trick  
\cite{fio}. Let us briefly recall the details.  
   
It suffices to show $\|\varrho(A)f \|^2 \le C( \|f\|^2 + \|\varrho(D)  
f\|^2)$, for some constant $C$ independent of $f$. Choose $C_1 > 0$  
with $|\sigma_m(A)|^2 \le C_1 |\sigma_m(D)|^2$. This is possible  
because $\sigma_m(A) \sigma_m(D)^{-1}$ is defined and continuous on  
the sphere bundle $S^{*}(\GR)$ of $A^*(\GR)$, a compact space.  Let $b  
> 0$ be smooth with $b^2 = (C_1 + 1) |\sigma_m(D)|^2 -  
|\sigma_m(A)|^2$ (this is defined only outside the zero section), and  
let $B \in \tPS{m}$ be an operator with principal symbol  
$\sigma_{m}(B)=b$. Then  
$$   
        (C_1 + 1) D^*D - A^*A - B^*B = R,   
$$   
with $R$ of order $l \le 2m - 1$. By replacing $B$ with $B_1$ such  
that $B_1 - B$ has order $l - m$ and $\sigma_{l-m}(B_1 - B) =  
\sigma_l(R)/2b$, we obtain that the order of the operator 
$(C_1+1)  
D^*D - A^*A - B_1^*B_1$ is less than $l$. Continuing in this way, we  
may assume that $R$ has order $\le 0$, so in particular $\varrho(R)$  
is bounded. Then  
$$   
        \|\varrho(A)f\|^2 \le (C_1+1) (\varrho(D^*D)f,f) - (\varrho(R)f,f) \le   
        C(\|\varrho(D)f\|^2 + \|f\|^2)   
$$   
for $C := \max\{ \|\varrho(R)\|, C_1 + 1\}$.   
\end{proof}

\begin{lemma}\label{lemma.5}\ Let $D =D^* \in \tPS{m}$, be elliptic    
with $\sigma_m(D) > 0$. Then we can find $C \ge 0$ such that, for any  
bounded representation $\varrho$ of $\tPS{\infty}$ on $\mathcal H$ we  
have $(\varrho(D) f,f) \ge -C(f,f)$, for all $f \in  
\cHu:=\varrho(\tPS{-\infty})\cH$.  
\end{lemma}

\begin{proof}\ The statement follows from the boundedness of $D$, if   
$m \le 0$, so assume that $m > 0$. Then the proof is the same as that  
of the previous lemma if in the proof of that lemma we replace $D^*D$  
with $D$ and take $A = 0$.  
\end{proof}

For the rest of the proof of Theorem \ref{theorem.red}, we shall fix a  
non-degenerate representation $\varrho$ of $\tPS{\infty}$ on  
${\mathcal H}$, and we shall identify the elements of $\tPS{\infty}$  
with unbounded operators with common domain $\cHu$, and the elements  
of $\alg\GR$ with bounded operators on $\mathcal H$.

\begin{corollary}\label{cor.1}\    
If $D = D^* \in \tPS{m}$, $m > 0$, is such that $\sigma_m(D) > 0$,  
then there exists $C \geq 0$ such that $D + C E \ge E$.  For any such $C  
\geq 0$ and any $A\in\tPS{m}$, the operator $A(D + C)^{-1}$ extends  
uniquely to a bounded operator on $\cH$.  
\end{corollary}

\begin{proof}\    
The first statement follows from the previous lemma. Fix $C \geq 0$ such   
that $D + CE \geq E$.  Lemma \ref{lemma.4} gives $\|Af\| \le C_1(\| f   
\| + \| (D + CE) f \|),$ for some $C_1 > 0$ and all $f \in \cHu$   
Consequently, there is $C_{2}>0$ with  
   
\begin{equation}\label{eq.ineq}   
        \|Af\| \le C_2 \| (D + C)f \|,   
\end{equation}   
for all $f\in\cHu$. Since $D+CE$ is essentially self-adjoint by  
Proposition \ref{ess} and $D+CE\geq E$, its range $\cH_{1}:=(D+CE)\cH$  
is dense by \cite[Theorem X.26]{resi2}. By \eqref{eq.ineq}, we obtain  
for $g=(D+CE)f\in\cH_{1}$  
$$   
        \| A(D + CE)^{-1}g \| \le C_2 \| g \|\,,   
$$    
which completes the proof.    
\end{proof}

\begin{corollary}\label{cor.2}\    
Consider now two self-adjoint, elliptic elements $D_1, D_2 \in  
\tPS{m}$, $m > 0$, with $D_i \ge E$ and $\sigma_m(D_i) > 0$, $i =  
1,2$.  Then $D_1D_2^{-1}$ extends uniquely to a bounded invertible  
operator.  
\end{corollary}

\begin{proof}\ By the previous corollary, both $D_1D_2^{-1}$ and $D_2D_1^{-1}$    
extend to bounded operators.  
\end{proof}

\begin{lemma}\label{lemma.6}\    
Let $D \in \tPS{m}$, $m >0$, be such that $D \ge E$, $\sigma_m(D) >  
0$, and $D^{-1} \in \ideal\GR$. Then we have $PD^{-k},D^{-k}P \in  
\alg\GR$, if $P$ has order $\le km$. Moreover, we have
$\sigma_{0}(PD^{-k})=\sigma_{km}(P)\sigma_{m}(D)^{-k}$, and
$\sigma_{0}(D^{-k}P)=\sigma_{m}(D)^{-k}\sigma_{km}(P)$.
\end{lemma}

\begin{proof}\ We notice that if $D$ satisfies the assumptions of the   
lemma, then $D^k$ satisfies them as well. We can assume then that   
$k = 1$.   
   
We shall check only that $PD^{-1} \in \alg \GR$. The relation $D^{-1}P  
\in \alg \GR$ can be proved in the same way or follows from the first  
one by taking adjoints.  
   
Let $A \in \tPS{-m}$ be with $AD -E = R \in \tPS{-\infty}$, $B_n \in  
\tPS{-\infty}$ be a sequence converging to $D^{-1}\in\ideal{\GR}$, and  
define $A_n := A - RB_n \in\tPS{-m}$. Then we have $A_n - D^{-1} =  
R(D^{-1} - B_n)$, thus, $PD^{-1}=PA_{n}-PR(D^{-1}-B_{n})$ first  
defined on the dense subspace $D\cHu$, has a unique bounded extension  
with $PD^{-1}\in\alg{\GR}$ because of $PA_n \in \tPS{0}$ and  
$$   
        \|PA_n - PD^{-1}\| \le \|PR\| \| D^{-1} - B_n\| \to 0\,, \quad  
        n \to \infty.  
$$   
Since $\sigma_{0}(PD^{-1})$ is the limit of 
$\sigma_{0}(PA_{n})=\sigma_{m}(P)\sigma_{-m}(A)=
\sigma_{m}(P)\sigma_{m}(D)^{-1}$, we obtain the formula for
the principal symbol as well.
\end{proof}

\begin{lemma}\label{lemma.7}\ Let $D \in \tPS{m}$ be with  $D \ge E$ and    
$\sigma_m(D) > 0$. Then $D^{-1} \in \ideal\GR$.  
\end{lemma}

\begin{proof}\    
{}From Lemma \ref{lemma.3} we know that $(E + \Delta)^{-m}$ is in  
$\ideal\GR$. By Corollary \ref{cor.2}, applied to $D_{1}=(E +  
\Delta)^{m}$ and $D_{2}=D^{2}$, we have $(E +  
\Delta)^{m}D^{-2}=D_{1}D_{2}^{-1}\in\alg\GR^{-1}$, thus, $D^{-2} = (E  
+ \Delta)^{-m}(E + \Delta)^{m}D^{-2} \in \ideal\GR$. Taking square  
roots completes the proof.  
\end{proof}

We now complete the proof of Theorem \ref{theorem.red}.

\begin{proof}\ Assume first that $m = 1$. Then the theorem follows from   
Lemma \ref{lemma.6} and Lemma \ref{lemma.7}. For arbitrary $m$,  
$D^{-1}\in \ideal\GR$, and hence we get $D^{-s} \in \ideal\GR$ by  
using functional calculus with continuous functions. A look at Lemma  
\ref{lemma.6} completes the proof.  
\end{proof}

We now obtain some corollaries of Theorem \ref{theorem.red}.  
   
For the following results, we need to define Sobolev spaces. Fix a  
metric on $A(\GR)$. Let then $\Delta:=\Delta_0 \in \Diff{\GR}$ be the  
Hodge-Laplacian acting on functions, and $\varrho$ be a  
non-degenerate, bounded representation as above. Then  
$D:=\varrho(E+\Delta)$ is essentially self-adjoint and strictly  
positive, hence we can define $D^s$, for each $s \in \RR$, using the  
functional calculus for essentially self-adjoint operators.  Then  
$H^s(\mathcal H,\varrho)$, {\em the $s$th Sobolev space of $(\mathcal  
H,\varrho)$}, is by definition, the domain of $D^{s/2}$ with the graph  
topology, if $s \ge 0$, or its dual if $s < 0$.

\begin{corollary}\label{cor.Sobolev}\    
The spaces $H^s(\mathcal H,\varrho)$ do not depend on the choice of  
the metric on $A(\GR)$, and every pseudodifferential 
operator $P \in \tPS{m}$ gives rise  
to a bounded map  
$H^{s}(\mathcal H ,\varrho) \to H^{s - m}(\mathcal H ,\varrho)$.  
\end{corollary}

\begin{proof}\ If we change the metric on the compact space $M$,   
we obtain a new Laplace operator, and $D$ will be replaced by a  
different operator $D_1$. However, by Corollary \ref{cor.2}, $D^s  
D_1^{-s}$ and $D_1^s D^{-s}$ are bounded for all even integer $s$. By  
interpolation, they are bounded for all $s$.  This proves the  
independence of the Sobolev space on the choice of a metric on $M$.  
   
The last claim follows from Lemma \ref{lemma.6} if $s$ is an integer.  
Let $H^{\infty} := \bigcap H^{k}(\mathcal H ,\varrho)$.  Then  
$\varrho(P)(H^{\infty})\subseteq H^{\infty}$.  Using this fact and  
applying the Phragmen-Lindel\"of principle to $s \mapsto  
(\varrho(D)^{s} P \varrho(D)^{-s} \xi,\xi')$, with $\xi,\xi' \in  
H^{\infty}$, we obtain the desired result for all $s$.  
\end{proof}

Similarly, we prove the following corollary. 
 
\begin{corollary} 
\label{coriso} 
Let $A\in\tPS{k}$ be elliptic. 
Then $\Lambda:=\varrho(E+A^{*}A)$ is essentially self-adjoint, and 
$\Lambda^{t}$ induces for all  $s,t\in\RR$  isomorphisms 
$$\Lambda^{t}:H^{s}(\cH,\varrho)\longrightarrow H^{s-2kt}(\cH,\varrho)\,.$$ 
\end{corollary}

Another corollary is related to the Cayley transform.

\begin{corollary}\label{cor.Cayley}\ If $A=A^{*} \in \tPS{m}$, $m >0$,    
is elliptic, then the Cayley transform $(A + iE)(A - iE)^{-1}$ 
of $A$ belongs to  $ \alg\GR$.  Moreover, we have
$$\sigma_{0}((A+iE)(A-iE)^{-1})=\sigma^{m}(A)\sigma_{m}(A)^{-1}=1\,,$$
where the last equality holds in the scalar case only.
\end{corollary}

\begin{proof}\    
We have $(A + iE)(A - iE)^{-1} = (A + iE)^2 (A^2 + E)^{-1} \in  
\alg\GR$, by Theorem \ref{theorem.red}, because $A^2 + E \ge E$ and  
$\sigma_{2m}(A^2 + E) > 0$.  The identity for the principal symbol follows from
the corresponding one in Lemma \ref{lemma.6}. 
\end{proof}

The Cayley transform of $A$ will be denoted in the following sections  
simply by $(A + iE)(A - iE)^{-1}$, because no more confusions can  
arise.

\section{Operators on open manifolds\label{Sec.O.M}}

One of the main motivations for studying algebras of  
pseudodifferential operators on groupoids is that they can be used to  
analyze geometric operators on certain complete Riemannian manifolds  
$(M_0,g)$ (without corners). The groupoids $\GR$ used to study these  
geometric operators will be of a particular kind.  They will have as  
space of units a compactification $M$ of $M_0$ to a manifold with  
corners such that $M_0$ will be an open invariant subset of $M$ with  
the property that the reduction of $\GR$ to $M_0$ is the product  
groupoid. If $M_0$ happens to be compact, then $M = M_0$, and our  
results simply reduce to the usual ``elliptic package'' for compact  
smooth manifolds without corners.  Our results thus can be viewed as a  
generalization of the classical elliptic theory from compact manifolds  
to certain non-compact, complete Riemannian manifolds.  
    
We now make explicit the hypothesis we need on the groupoid $\GR$.  
\vspace*{2mm}

{\bf Assumptions.}\ {\em In this and the following sections, $M_0$  
will be a smooth manifold without corners which is diffeomorphic to  
(and will be identified with) an open dense subset of a compact  
manifold with corners $M$, and $\GR$ will be a differentiable groupoid  
with units $M$, such that $M_0$ is an invariant subset and}  
$$   
        \GR_{M_0}\cong M_0 \times M_0.   
$$   
     
The above assumptions have a number of useful consequences for $\GR$,  
$M$, and $M_0$, and we shall use them in what follows, without further  
comment.  
   
Let $A = A(\GR)$. First of all, $A\vert_{M_0} \cong TM_0$. Fix a  
metric on $A$.  The metric on $A$ then restricts to a metric on $M_0$,  
so $M_0$ is naturally a Riemannian manifold such that the map $r:\GR_x  
\to M_0$ is an isometry for any $x \in M_0$. Moreover, because $M$ is  
compact, all metrics on $M_0$ obtained by this procedure will be  
equivalent:\ if $g_1$ and $g_2$ are metrics on $M_0$ obtained from  
metrics on $A$, then we can find $C,c > 0$ such that $cg_1 \le g_2 \le  
Cg_1$ (this is of course not true for any two metrics on the  
non-compact smooth manifold $M_0$).  The same result holds true for  
the induced smooth densities (or measures) on $M_0$, and hence all the  
spaces $L^2(M_0)$ defined by these measures actually coincide.  
    
Let $\pi$ be the vector representation of $\tPS{\infty}$ on $\CI(M)$  
(uniquely determined by $(\pi(P)f) \circ r = P (f \circ r)$, see  
Equation \eqref{eq.vector.rep}). Then $\pi(\tPS{\infty})$ maps  
$\CIc(M_0)$ to itself.  Fix $x \in M_0$. The regular representation  
$\pi_x : \tPS{\infty} \to End(C_c^{\infty}(\GR_x))$ is equivalent to  
$\pi$ via the isometry $r: \GR_x \to M_0$, and hence $\pi$ is a  
bounded representation of $\tPS{\infty}$ on $L^2(M_0)$.  
   
We now relate the geometric operators on $M_0$, defined using a metric  
induced from $A$, and the geometric operators on $\GR$, $\GR$ as above  
($\GR_{M_0} \cong M_0 \times M_0$).  
   
We start with a $\Cliff(A)$-module $W$ on $M$ together with an  
admissible connection $\nabla^W \in \Diff{\GR;W,W\otimes A^*}$,  
defined as an invariant family of differential operators on $\GR_x =  
d^{-1}(x)$.  Fix $x \in M_0$ arbitrary. Then the restriction of  
$r^*(W)$ to $\GR_x$ is a Clifford module on $\GR_x$, which hence can  
be identified with a Clifford module $W_0$ on $M_0$, using the  
isometry $\GR_x \cong M_0$.  
   
Let $\Dir^{W} \in \Psi^{1}(\GR;W)$ be the Dirac operator on $\GR$  
associated to $W$ and its admissible connection, and let $\Dir^{W_0}$  
be the Dirac operator on $M_0$ associated to $W_0$ and its admissible  
connection obtained by pulling back the connection on $\GR_x$. These  
operators are related as follows.

\begin{theorem}\label{Theorem.Image}\ The Dirac operator $\Dir^W$ on $\GR$ acts   
in the vector representation as $\Dir^{W_0}$, the Dirac operator on  
$M_0 \subset M$ defined above. More precisely,  
$$   
        \pi(\Dir^{W}) = \Dir^{W_0}.  
$$   
\end{theorem}

\begin{proof}\    
By construction, $\Dir^{W_0}$ is, up to similarity, the restriction of  
$\Dir^W$ to one of the fibers $\GR_x$, with $x \in M_0$.  
\end{proof}

At first sight, the above theorem applies only to a very limited class  
of (admissible) Dirac operators on $M_0$, the ones coming from  
$\Cliff(A)$-modules. Not every Dirac operator on a Clifford module on  
$M_0$ can be obtained in this way. However, as we shall see in a  
moment, if we are given a Clifford module on $M_0$, we can always  
adjust our compatible connection so that the resulting Dirac operator  
comes from a Dirac operator on $\GR$ (corresponding to a  
$\Cliff(A)$-module).

\begin{theorem}\label{Theorem.Geometric}\    
Suppose there exists a compact subset $M_1 \subset M_0$ which is a  
deformation retract of $M$. Let $W_0$ be a Clifford module on  
$M_0$. Then we can find an admissible connection on $W_0$ such that  
the associated admissible Dirac operator $\Dir^{W_0}$ is (conjugate  
to) $\pi(\Dir^{W})$, for some $\Cliff(A)$-module $W$, $\Dir^W$ being  
the Dirac operator on $\GR$ associated to $W$.  
   
If $W_0$ is a spin bundle, then we can choose this connection to   
be the Levi-Civita connection on $W_0$.   
\end{theorem}

\begin{proof}\ Using the deformation retract $f : M \to M_1$, we define   
(up to isomorphism) $W = f^*(W_0)$. Then $W\vert_{M_0} \simeq W_0$,  
the isomorphism being uniquely determined up to homotopy. Moreover, we  
have a (non-canonical) isomorphism $A \simeq f^*(TM_0)$ of vector  
bundles, which allows us to define a $\Cliff(A)$-module structure on  
$W$. By replacing $W_0$ with an isomorphic bundle, we can assume then  
that $W_0 = W\vert_{M_0}$, as Clifford modules. Choose an admissible  
connection on $W$. Theorem \ref{Theorem.Image} then gives that  
$\pi(\Dir^{W}) = \Dir^{W_0}$.  
\end{proof}

\section{Spectral properties\label{Sec.SP}}

We shall use now the results of the previous section to study  
operators on suitable Riemannian manifolds.  We are interested in  
spectral properties, Fredholmness, and compactness for these  
operators.  The results of this section extend essentially without  
change to the case of families of such manifolds.  
   
{\em We fix, throughout this section, a groupoid $\GR$ satisfying the  
assumptions of Section \ref{Sec.O.M}.} In particular, $M_0$ is an open  
invariant subset of $M$ and $\GR_{M_0}\cong M_0 \times M_0$.  
   
We denote as before by $\alg\GR$ the closure of $\tPS0$ in the norm  
$\|\;\cdot\|$ and by $\ideal\GR$ the closure of $\tPS{-\infty}$ in the  
same norm.  Our analysis of geometric operators on $M_0$ depends on  
the structure of the algebras $\alg\GR$ and $\ideal\GR$.  The results  
of Section \ref{Sec.BR} applied to our groupoid $\GR$ (satisfying the  
assumptions of Section \ref{Sec.O.M}) give the following. Let $\mI =  
\ideal{\GR_{M_0}}$, then $\mI$ is isomorphic to  
$\mathcal{K}(L^2(M_0))$, the algebra of compact operators on $L^2(M_0)  
= L^2(M)$, the isomorphism being induced by the vector representation  
$\pi$, or by any of the representations $\pi_x$, $x \in M_0$, and the  
isometry $\GR_x \simeq M$. Otherwise, if $x \notin M_0$, then $\pi_x$  
descends to a representation of $Q(\GR) : = \alg\GR/\mI$.  
   
We shall study various spectra, for this purpose, the results of  
Section \ref{Sec.Sobolev} will prove indispensable.  
   
We denote by $\sigma(P)$ the spectrum of an element $P \in \alg\GR$  
and by $\sigma_{Q(\GR)}(P)$ the spectrum of the image of $P$ in  
$Q(\GR) := \alg\GR/\mI = \alg\GR/\ideal{\GR_{M_0}}$.  These  
definitions extend to elliptic, self-adjoint elements $P \in \tPS{m}$,  
$m >0$, using the Cayley transform, as follows. Let $f (t) = (t +  
i)/(t -i)$ and $f(P):=(P + i)(P - i)^{-1} \in \alg{\GR}$ be its Cayley  
transform, which is defined by Corollary \ref{cor.Cayley}.  We define  
then  
\begin{equation}   
        \sigma(P) := f^{-1}(\sigma(f(P)))\,, \quad  
        \mbox{ and } \quad   
        \sigma_{Q(\GR)}(P) := f^{-1}(\sigma_{Q(\GR)}(f(P))).   
\end{equation}   
We observe that if $P$ is identified with an unbounded, self-adjoint  
operator on a Hilbert space, then the relation $\sigma(P) :=  
f^{-1}(\sigma(f(P)))$ is automatically satisfied, by the spectral  
mapping theorem.  
    
The spectrum and essential spectrum of an element $T$ acting as an  
unbounded operator on a Hilbert space will be denoted by $\sigma(T)$  
and, respectively, by $\sigma_{ess}(T)$.  (We shall do that for an  
operator of the form $T = \pi(P)$, with $\pi$ the vector  
representation and $P \in \PSF{m}$, $m > 0$, elliptic, self-adjoint.)  
   
We shall formulate all results below for operators acting on vector  
bundles. Fix an elliptic operator $A \in \PSF{m}$, $m >0$, then for any  
$P \in \PSF{k}$, we have $P_1 := P (E + A^*A)^{-k/2m} \in \alg{\GR;F}$,  
by Theorem \ref{theorem.red}.  
    
Let us notice that if $\pi$ is the vector representation of $\tPS{0}$  
on $L^2(M) = L^2(M_0)$ (see Equation \ref{eq.vector.rep} for the  
definition of the vector representation), then the spaces $H^s(M) =  
H^s(L^2(M),\pi)$ are the usual Sobolev spaces associated to the  
Riemannian manifold of bounded geometry $\GR_x \simeq M_0$  
\cite{roe1,shubin92}.  If we are working with sections of a Hermitian  
vector bundle $F$, then we write $H^s(M;F) := H^s(L^2(M;F),\pi)$.

\begin{theorem}\label{theorem.I}\    
Let $M_0 \subset M$, the groupoid $\GR$, and $A \in \PSF{k}$, elliptic,  
be as above.  
   
(i)\ If $P \in \PSF{m}$ is such that the image of $P_1:=P(E +  
A^*A)^{-m/2k} \in \alg\GR$ in $Q(\GR) : = \alg\GR/ \ideal{\GR_{M_0}}$  
is invertible, then $\pi(P)$ extends to a Fredholm operator  
$H^{m}(M;F) \to L^2(M;F)$.  
   
(ii)\ If $P_{1}:=P(E + A^*A)^{-m/2k}$ maps to zero in $Q(\GR)$, then  
$\pi(P)$ is a compact operator $H^{m}(M;F) \to L^2(M;F)$.  
   
(iii)\ If $P \in \PSF{0}$ or $P \in \PSF{m}$, $m > 0$, is self-adjoint,  
elliptic, then $\sigma(\pi(P)) \subseteq \sigma(P)$ and  
$\sigma_{ess}(\pi(P)) \subseteq \sigma_{Q(\GR)}(P)$.  
\end{theorem}

\begin{proof}\ Let $P_1 := P(E + A^*A)^{-m/2k}$, as above.    
   
(i)\ Choose $Q_1 \in \alg\GR$ such that  
$$  
        Q_1P_1 - E, P_1 Q_1 - E \in \mI := \ideal{\GR_{M_0}},  
$$   
and define $Q = \pi ((E + A^*A)^{-m/2k}Q_1)$. Then $ \pi(P) Q -  
\id_{L^{2}(M)} \in \pi(\mI) = \mathcal K$. Similarly, we find a right  
inverse for $\pi(P)$ up to compact operators. Thus, $\pi(P)$ is  
Fredholm.  
   
(ii)\ The operator $\pi(P) : H^{m}(M;F) \to L^2(M;F)$ is the product  
of the bounded operator $\pi(E + A^*A)^{m/2k} : H^{m}(M;F) \to  
L^2(M;F)$ and of the compact operator $\pi(P_1)$.  
    
(iii)\ For $P$ in a $C^*$-algebra $A_0$ and $\varrho$ a bounded 
$*$-representation of $A_0$, the spectrum of $\varrho(P)$ is contained 
in the spectrum of $P$ (we do not exclude the case where they are 
equal). If $P \in \PSF0$, this gives (iii), by taking $A_0 = \alg\GR$ 
or $A_0 = Q(\GR)$. If $P \in \PSF{m}$, $m >0$, is self-adjoint, 
elliptic, then we use the result we have just proved for $f(P) =(P + 
iE) (P - iE)^{-1}$, the Cayley transform of $P$. 
\end{proof}

It is interesting to observe the following. Both (i) and (ii) can be  
proved using (iii). However, because (i) and (ii) are more likely to  
be used, we also included separate, simpler proofs of (i) and (ii). A  
derivation of (i) and (ii) from (iii) can be obtained as in the proof  
of the following theorem.  
  
Let $\pi:\tPS{\infty}\longrightarrow \End(\cunc(M))$ be the vector  
representation. Then the homogeneous principal symbol $\sigma_{0}(P)$  
of $P\in\tPS{0}$ can be recovered from the action of $\pi(P)$ on  
$\CI(M)$ by oscillatory testing as in the classical case.  Indeed, let  
$x\in M$ and $\xi\in A_{x}^{*}(\GR)=T^{*}_{x}\GR_{x}$ be  
arbitrary. Then we have  
\begin{equation}\label{psy}  
        \sigma_{0}(P)(\xi)= \lim_{t\to\infty}\big[  
        e^{-itf}\pi(P)\varphi e^{itf}\big](x)  
\end{equation}  
for all $\varphi\in\cunc(M)$ with $\varphi=1$ near $x$, and all 
$f\in\CI(M,\RR)$ with $d(f\circ r)\neq 0$ on $\supp(\varphi\circ r)$, 
and $d(f\circ r)|_{x}=\xi$.  The proof of \eqref{psy} uses the fact 
that the homogeneous principal symbol 
($\sigma_{0}(P)(\xi)=\sigma_{0}(P_{x})(\xi)$) as well as the action of 
$\pi(P)$ on $\cunc(M)$ ($(\pi(P)h)(x)=P_{x}(h\circ r)|_{\GR_{x}}(x)$) 
are defined using the manifold $\GR_{x}$ only. Thus, the classical 
result applies. 
 
Suppose now that the vector representation $\pi : \ideal\GR \to  
\mathcal{B}(L^2(M))$ is injective. Then the above result can be  
sharpened to a necessary and sufficient condition for Fredholmness, 
respectively for compactness.  We first note that since the principal 
symbol of a pseudodifferential operator can be determined from its 
action on functions, the representation $\pi : \alg\GR \to 
\mathcal{B}(L^2(M))$ is also injective. Indeed, this follows from  
formula \eqref{psy}.

\begin{theorem}\label{theorem.II}\    
Assume that the vector representation $\pi$ is injective on $\ideal\GR$.  
Using the notation from the above theorem, we have.  
   
(i)\ If $P \in \PSF{m}$ is such that $\pi(P)$ defines a Fredholm  
operator $H^{m}(M;F) \to L^2(M;F)$ then the image of $P(E +  
A^*A)^{-m/2k} \in \alg\GR$ in $Q(\GR) : = \alg\GR/ \ideal\GR$ is  
invertible.  
   
(ii)\ If $\pi(P)$ defines a compact operator $H^{m}(M;F) \to  
L^2(M;F)$, then the image of $P(E + A^*A)^{-m/2k}$ in $Q(\GR)$  
vanishes.  
   
(iii)\ If $P \in \PSF{0}$ or $P \in \PSF{m}$, $m > 0$, is self-adjoint,  
elliptic, then $\sigma(\pi(P)) =\sigma(P)$ and $\sigma_{ess}(\pi(P)) =  
\sigma_{Q(\GR)}(P)$.  
\end{theorem}

\begin{proof}\  An injective representation $\pi$ of   
$C^*$-algebras preserves the spectrum, and in particular, $a$ is  
invertible if, and only if, $\pi(a)$ is invertible.  
   
Denote by $\mathcal B$ the algebra of bounded operators on $L^2(M;F)$.  
The morphism $\pi':Q(\GR) \to \mathcal{B/K}$ induced by $\pi$ is also  
injective.  Fix $P_0 \in \alg\GR$.  Then $P_0$ is invertible if, and  
only if, $\pi(P_0)$ is invertible. By replacing $P_0$ with $P_0 -  
\lambda E$, we obtain $\sigma(P_0) = \sigma(\pi(P_0))$.  We see then  
that $P_0$ is invertible modulo $\ideal{\GR_{M_0}}$ if, and only if,  
$\pi(P_0)$ is invertible modulo compact operators.  This gives  
$\sigma_{Q(\GR)}(P_0) = \sigma_{ess}(\pi(P_0))$.  We thus obtain (iii)  
if we take $P_0 = P$ or $P_0 = f(P)$, the Cayley transform of $P$.  
   
(i) By definitions, $\pi(P_1)$ is Fredholm if, and only if, $\pi(P)$  
defines a Fredholm operator $H^{m}(M;F) \to L^2(M;F)$. Then  
\begin{eqnarray*}   
        \pi(P_1) \text{ is Fredholm } &\iff &0  
        \notin\sigma_{ess}(\pi(P_1)) \\ &\iff& 0 \notin  
        \sigma_{Q(\GR)}(P_1) \\ &\iff& P_1 \text{ is invertible in }  
        Q(\GR)\,.  
\end{eqnarray*}   
   
For (ii), a similar reasoning holds:  
\begin{eqnarray*}   
        \pi(P_1) \text{ is compact }& \iff &  
        \sigma_{ess}(\pi(P_1^*P_1)) = \{0\} \\ &\iff &  
        \sigma_{Q(\GR)}(P_1^{*}P_{1}) = \{0\} \\ &\iff & P_1 = 0 \in  
        Q(\GR)\,.  
\end{eqnarray*}   
\end{proof}

The criteria in the above theorems can be made even more explicit in   
particular examples.

\begin{theorem}\label{theorem.III}\   
Suppose the restriction of $\GR$ to $M \smallsetminus M_0$ is amenable,  
and the vector $\pi$ representation is injective. Then,  
   
(i)\ $P : H^s(M;F) \to L^2(M;F)$ is Fredholm if, and only if, $P$ is  
elliptic and $\pi_x(P): H^s(\GR_x,r^*F) \to L^2(\GR_x,r^*F)$ is  
invertible, for any $x \not \in M_0$.

(ii)\ $P : H^s(M;F) \to L^2(M;F)$ is compact if, and only if, its  
principal symbol vanishes, and $\pi_x(P) = 0$, for all $x \not \in  
M_0$.  
   
(iii)\ For $P \in \PSF{0}$, we have  
      $$\sigma_{ess}(\pi(P)) = \bigcup_{x  
     \not \in M_0}\sigma(\pi_x(P)) \cup  
      \bigcup_{\xi\in S^{*}\GR}{\rm spec}(\sigma_{0}(P)(\xi))\,,$$ 
      where ${\rm spec} (\sigma_{0}(P)(\xi))$ denotes the spectrum 
      of the linear map $\sigma_{0}(P)(\xi):E_x\rightarrow E_x$. 
   
(iv)\ If $P \in \PSF{m}$, $m >0$, is self-adjoint, elliptic, then we  
have $\sigma_{ess}(\pi(P)) = \cup_{x \not \in M_0}\sigma(\pi_x(P))$.  
\end{theorem}

\begin{proof}\ Again, (i) and (ii) follow from (iii) and (iv). The assumption   
$\alg\GR = \ralg\GR$ implies $\alg\GR/\mI = \ralg\GR/\mI$. Because the  
groupoid obtained by reducing $\GR$ to $M \smallsetminus M_0$ is  
amenable, the representation $\varrho := \prod \pi_x$, $x \not \in  
M_0$ is injective on $Q(\GR)$. This gives $\sigma_{Q(\GR)}(T) = \cup_x  
\sigma(\pi_x(T))$, $x \not \in M_0$, for all $T \in \alg\GR$.  
\end{proof}

Another explicit criterion is contained in the theorem below.

\begin{theorem}\label{theorem.IV}\    
Suppose the vector representation $\pi$ is injective and $M \setminus  
M_0$ can be written as a union $\bigcup_{j=1}^rZ_j$ of closed,  
invariant manifolds with corners $Z_j \subset M$.  
   
(i)\ Let $P\in \PSF{m}$, then $P : H^s(M) \to L^2(M)$ is Fredholm if,  
and only if, it is elliptic and $\inn_{Z_j}(P) : H^s(Z_j) \to  
L^2(Z_j)$ is invertible, for all $j$.

(ii)\ Let $P\in \PSF{m}$, then $P : H^s(M) \to L^2(M)$ is compact if,  
and only if, its principal symbol vanishes and $\inn_{Z_j}(P) = 0$,  
for all $j$.  
   
(iii)\ For $P \in \PSF{0}$, we have $$\sigma_{ess}(\pi(P)) =  
       \bigcup_{j=1}^r\sigma(\inn_{Z_j}(P)) \cup   
      \bigcup_{\xi\in S^{*}\GR}{\rm spec}(\sigma_{0}(P)(\xi))\,.$$ 
          
(iv)\ Suppose $P \in \PSF{m}$, $m >0$, is self-adjoint, elliptic. Then  
we have $\sigma_{ess}(\pi(P)) = \cup_{j=1}^r\sigma(\inn_{Z_j}(P))$.  
\end{theorem}

\begin{proof}\ The representation 
$\alg\GR/\mI \to \oplus_j \alg{\GR_{Z_j}}\oplus\mathcal{C}(S^{*}\GR;\End(F))$ 
given by 
the restrictions $\inn_{Z_{j}}$ and the homogeneous principal symbol   
is injective. This gives (iii) and (iv).  
For  $m>0$ note that we have $\sigma_{0}(f(P))=\id_{F}$ for  the Cayley
transform $f(P)=(P+iE)(P-iE)^{-1}\in\alg{\GR}$ of $P$, 
and $f^{-1}(1)=\{\infty\}$.

To obtain (i) and (ii) from  
(iii) as above, it is enough to observe that the operator $P_1 = P(E +  
A^*A)^{-m/2k}$, (with $A$ elliptic of order $m$, fixed) belongs to $\mI =  
\ideal{\GR_{M_0}}$ if, and only if, $\inn_{Z_j}(P_1) = 0$ for all $j$,  
and that $\inn_{Z_j}(P_1) = 0$ if, and only if, $\inn_{Z_j}(P) = 0$.  
\end{proof}

In Section \ref{Sec.Examples.III}, we shall see  examples of  
groupoids for which the conditions of the above theorem are  
satisfied. Having this natural characterization of Fredholmness, it is  
natural to ask for an index formula for these operators, at least in  
the case when the restriction of $\GR$ to each component $S$ of is  
such that $r(\GR_x)$ has constant dimension for all $x$ in a fixed  
component $S$ of $Y_k \smallsetminus Y_{k-1}$ (that is, when the  
restriction of $A(\GR)$ to each component $S$ of $Y_k \smallsetminus  
Y_{k-1}$ is a regular Lie algebroid). The results of  
\cite{NistorIndFam} deal with a particular case of this problem, when  
$M = Y_n$, the induced foliation on $M$ is a fiber bundle, and the  
isotropy bundle can be integrated to a bundle of Lie groups that  
consists either of compact, connected Lie groups or of  
simply-connected, solvable Lie groups.

\section{Examples III: Applications\label{Sec.Examples.III}}

For geometric operators $P$, the operators $\pi_x(P)$ and $\pi_Z(P)$  
appearing in the statements of the above theorems are again geometric  
operators of the same kind (Dirac, Laplace, ...). This leads to very  
explicit criteria for their Fredholmness and to the inductive  
determination of their spectrum.  
   

\begin{example}\ If $\GR = M \times M$ is the pair groupoid, then     
$\ideal\GR \cong \mathcal K = \mathcal K(L^2(M))$ and all the results stated    
above were known for these algebras. In particular, the exact sequence   
$$   
        0 \to \mathcal K \to \alg\GR \to \mathcal{C}(S^*M) \to 0   
$$   
is well-known.  Moreover, the criteria for compactness and  
Fredholmness are part of the classical elliptic theory on compact  
manifolds. There is no need for an inductive determination of the  
spectrum in this case.  
\end{example}

We are now going to apply the results of the previous section to the  
$c_{n}$-calculus considered in Example \ref{excn}. The main result  
being an inductive method for the determination of the essential  
spectrum of Hodge-Laplace operators. Because the $b$-calculus  
corresponds to the special case $c_{H}=1$ for all boundary hyperfaces  
$H$ of $M$, we in particular answer an question of Melrose on the  
essential spectrum of the $b$-Laplacian on a compact manifold $M$ with  
corners \cite[Conjecture 7.1]{MelroseScattering}.  
   
Let $\GR(M,c)$ be the groupoid constructed in Example $\ref{excn}$ for  
an arbitrary system $c=(c_{H})$.

\begin{lemma}  
\label{lemma7}  
The groupoid $\GR(M,c)$ is amenable and the   
vector representation of $\alg{\GR(M,c)}$ is injective.   
\end{lemma}

\begin{proof}\ It is enough to prove that the representation   
$\pi$ is injective on $\ideal{\GR(M,c)}$, because we can recover the  
principal symbol of a pseudodifferential operator from its action on  
functions, as explained in the previous section.  
   
The groupoid $\GR$ is amenable because the composition series of  
Theorem \ref{Theorem.CS} are associated to the groupoids $(S \times S)  
\times \RR^k$, which are amenable groupoids.  
   
It is then enough to prove that each representation of the form  
$\pi_x$ is contained in the vector representation. Let $x \in F$ be an  
interior point. By considering a small open subset of $x$, we can  
reduce the problem to the case when the manifold $M$ is of the form  
$[0,\infty)^k \times \RR^{n-k}$. Then the result is reduced to the  
case $k=n=1$ using Proposition \ref{prop.tens}. But for this case  
$\ideal\GR$ is isomorphic to the crossed product algebra $C_0(\RR \cup  
\{-\infty\}) \rtimes \RR$ and the vector representation corresponds to  
the natural representation on $L^2(\RR)$ in which $C_0(\RR \cup  
\{-\infty\})$ acts by multiplication and $\RR$ acts by  
translation. This representation is injective (it is actually often  
used to define this crossed product algebra). {}From this the result  
follows.  
\end{proof}

The algebra $C_0(\RR \cup \{-\infty\}) \rtimes \RR$ is usually called  
the algebra of {\em Wiener-Hopf operators} on $\RR$, for which it is  
well-known that the vector representation is injective.  
   
Fix now a metric $h$ on $A=A(\GR(M,c))$, and let  
$\Delta_p^{c}:=\Delta_{p}^{\GR(M,c)}$ be the corresponding  
Hodge-Laplacian acting on $p$-forms. Note that each boundary hyperface  
$H$ of $M$ is a closed, invariant submanifold with corners, whereas  
the interior $M_{0}:=M\setminus\partial M = M\setminus\bigcup H$ is  
invariant and satisfies $\GR_{M_0}\cong M_{0}\times M_{0}$. We are in  
position then to use the results of the previous section.  
   
First we need some notation.  For each hyperface $H$ of $M$, we  
consider the system $c^{(H)}$ determined by $c^{(H)}_{F}=c_{F'}$ for  
all boundary hyperfaces $F'$ of $M$ with $F:=H\cap F'\neq \emptyset$,  
as in the Example \ref{excn}. By the construction of the groupoid  
$\GR(M,c)$,  
$$   
        \GR(M,c)_H \cong \GR(H,c^{(H)}) \times \RR.   
$$   
It will be convenient to use the Fourier transform to switch to the  
dual representation in the $\RR$ variable, so that the action of the  
group by translation becomes an action by multiplication. Then  
pseudodifferential operators on $\GR(M,c)_H$ become families of  
pseudodifferential operators on $\GR(H,c^{(H)})$ parametrized by  
$\RR$. Using also \eqref{lapl}, this reasoning then gives  
\begin{multline}\label{ide}   
        \inn_{H}(\Delta_{p}^{c})=\Delta_{p}^{\GR(M,c)_{H}} =  
        \begin{cases} \lambda^{2} + \Delta_{0}^{ c^{(H)} }, & \text{  
        if } p = 0, \\ \big (\lambda^{2} + \Delta_{p}^{ c^{(H)} } \big  
        )\oplus \big ( \lambda^{2} + \Delta_{p-1}^{ c^{(H)} } \big ),  
        & \text{ if } p > 0, \end{cases}  
\end{multline}   
because for $p > 0$ the space of $p$-forms on the product with  
$[0,\infty)$ splits into the product of the spaces of $p-1$ and $p$  
forms that contain, respectively, do not contain, $dt$, $t \in  
[0,\infty)$.  
   
Denote by $m_H^{(p)} = \min \sigma(\Delta_{p}^{c^{(H)}})$ and by  
$m^{(p)} = \min_H m_H^{(p)}$.  Then $m^{(p)}\geq0$ because the  
Hodge-Laplace operators $\Delta_{p}^{c^{(H)}}$ are positive operators.  
  
On the other hand, note that $\pi(\Delta_p^c)$ is (conjugated to)  
$\Delta_p$, the Hodge-Laplace operator acting on $p$-forms on the  
complete manifold $M_0 := M \smallsetminus \pa M$, with the induced  
metric from $A(M,c)$.

\begin{theorem}\label{mainth}\   
Consider the manifold $M_0$, which is the interior of a
compact manifold with corners $M$, with the metric induced from
$A(M,c)$. Then the essential spectrum of the (closure of the)
Hodge-Laplacian $\Delta_{p}$ acting on $p$-forms on $M_0$ is
$[m,\infty)$, with $m = m^{(0)}$, if $p = 0$, or $m = \min \{ m^{(p)}
, m^{(p-1)}\}$, if $p>0$, using the notation explained above.
\end{theorem}

In particular, the spectrum of $\Delta_{p}$ itself is the union of  
$[m,\infty)$ and a discrete set consisting of eigenvalues of finite  
multiplicity.  
   
\begin{proof}   
We are going to apply Theorem \ref{theorem.IV} (iv), with $Z_j$  
ranging through the set of hyperfaces of $M$; this is possible because  
of Lemma \ref{lemma7}. Furthermore, note that by the definition of the  
groupoid structure on $\GR(M,c)$ in Example \ref{excn}, the boundary  
hyperfaces $H$ of $M$ are closed, invariant submanifolds with  
$M\setminus M_{0}=\bigcup_{H}H$.  
  
For each boundary hyperface $H$ of $M$, we have by \eqref{ide}   
$$   
        \sigma(\inn_{H}(\Delta_{p}^{c})) =  
        \bigcup_{\lambda\in\RR}\left(
        \lambda^{2}+\left(\sigma(\Delta_{p}^{c^{(H)}})  
        \cup\sigma(\Delta_{p-1}^{c^{(H)}})\right)\right) =  
        [\min\{m_{H}^{(p-1)},m_{H}^{(p)}\},\infty)\,,  
$$   
where $\lambda^{2}+\sigma(\Delta_{p-1}^{c^{(H)}})$ is missing if  
$p=0$.  Since $\Delta_{p}$ is essentially self-adjoint and elliptic,  
Theorem \ref{theorem.IV} (iv) completes the proof.  
\end{proof}

Using an obvious inductive procedure, we then obtain the following
more precise result on the spectrum of the Laplace operator acting
on functions. 
\footnote{This refinement of our theorem was suggested to us by a
comment of Richard Melrose during the talk of the first named author
at the Oberwolfach meeting {\em Geometric Analysis and Singular
Spaces}, June 2000.}

\begin{corollary}\ Let $M_0$ be as above, then the spectrum of
the (closure of the) Laplace operator $\Delta_0$ on $M_0$ is
$\sigma(\Delta_0) = [0,\infty)$, and hence it coincides with its
essential spectrum.
\end{corollary}

\begin{proof}\ Let $F $ be a minimal face of $M$ (that is, not
containing any other face of $M$). Then $F$ is a compact manifold
without corners and hence the Laplace operator on $F$ contains $0$ in
its spectrum. The above theorem then shows that $[0,\infty) \subset
\sigma_{ess}(\Delta_0)$. On the other hand, $\Delta_0$ is positive,
and hence $\sigma(\Delta_0) \subset [0,\infty)$. This completes the
proof.
\end{proof}


\end{document}